\documentclass[graybox]{svmult}

\usepackage{mathptmx}       
\usepackage{helvet}         
\usepackage{courier}        
\usepackage{type1cm}        
\usepackage{makeidx}         
\usepackage{graphicx}        
\usepackage{multicol}        
\usepackage[bottom]{footmisc}
\usepackage{amssymb}
\usepackage{amsmath}

\newcommand{\ds}{\displaystyle}
\newcommand{\pder}[2]{ {\frac{\partial {#1}}{\partial {#2}}}}
\def\defeq{\mathrel{\mathop:}=}

\makeindex             

\begin{document}

\title*{An implicit boundary integral method for interfaces evolving by Mullins-Sekerka dynamics}
\titlerunning{Implicit boundary integral method for Mullins-Sekerka interface dynamics}
\author{Chieh Chen, Catherine Kublik and Richard Tsai}
\institute{Chieh Chen \at Goldman Sachs, 6011 Connection Drive, Irving TX, 75039 USA \email{Chieh.Chen@gs.com}
\and Catherine Kublik \at Department of Mathematics, University of Dayton,  300 College Park, Dayton OH, 45469 USA \email{ckublik1@udayton.edu} \and
Richard Tsai \at Department of Mathematics and Institute for Computational Engineering and Sciences, University of Texas at Austin, 2515 Speedway, 78712 Austin, USA. \\Department of Mathematics, KTH Royal Institute of Technology, SE-100 44, Stockholm, Sweden \email{tsai@kth.se} }
%
%
\maketitle

\abstract*{We present an algorithm for computing nonlinear interface dynamics driven by bulk diffusion using integral equations for implicitly defined interfaces. The computation of the dynamics involves solving Laplace's equation with Dirichlet boundary conditions on multiply connected and unbounded domains and propagating the interface using a normal velocity obtained from the solution of the PDE at each time step. Our method is based on a simple formulation which rewrites boundary integrals as volume integrals over the entire space. The resulting algorithm thus inherits the benefits of both level set method and boundary integral methods to simulate the nonlocal front propagation problem with possible topological changes. We present numerical results in both two and three dimensions to demonstrate the effectiveness of the algorithm.}


\section{Introduction}
\label{intro}

This paper proposes an algorithm for simulating interfacial motions in Mullins-Sekerka dynamics  using a boundary integral equations formulated on implicit surfaces. The proposed method is able to simulate the dynamics on unbounded domains, and as such will provide a tool to better understand the dynamics obtained in these situations. 
As pointed out in Gurtin's paper \cite{gurtin86},  the behavior of the dynamics are far more interesting on domains that are unbounded. However, it seems that there is no other computational method that can simulate the nonlinear interfacial motion of Mullins-Sekerka dynamics in three dimensions on unbounded domains, allowing the interfaces to merge, break up, without the need of finding explicit representations. 

The Mullins-Sekerka flow is a Stefan-type free boundary problem involving an nonlocal interfacial motion dynamically controlled by the solution of Laplace's equation with appropriate boundary conditions obtained on both sides of the interface. As a result, it is reasonable to consider boundary integral methods, especially for exterior domains, combined with level set methods \cite{osher_sethian88}\cite{Giga:2006}, for easy tracking and for being able to handle topological changes. 
%

Various numerical methods have been proposed to solve elliptic problems on irregular domains. We mention here some of the few "usual suspects" that use finite element methods  \cite{babuska70,dolbow_harari09,hansbo_hansbo02,hansbo_hansbo04,huang_zou02}, finite difference techniques \cite{bedrossian_teran10,chern_shu07, gibou_fedkiw05,johansen97,johansen_colella98,leveque_li94}, 
and boundary integral methods \cite{atkinson_chandler90, zhu_chen_hou96}. 
Finite element methods require an explicit representation (e.g. triangulation) of the domain which makes them less tractable for an evolving domain. 
There is also a wide range of finite difference methods for solving elliptic equations, some of them using the level set method \cite{gibou_fedkiw_cheng_kang02,liu_fedkiw_kang00, gibou_fedkiw05}. 
%
%
The Immersed Interface Method (IIM) \cite{leveque_li94,li_ito06}  is a popular technique for solving elliptic equations, particularly if the coefficients in the equation are discontinuous. This technique uses an adaptive finite difference scheme with a locally adaptive stencil. 
Unlike finite element method, the IIM can be used with an implicit representation of the interface. 
Nevertheless, both finite element and finite difference based techniques discretize the whole domain and therefore cannot be applied directly to "exterior" problems without imposing additional artificial boundary conditions. 

In contrast, boundary integral methods (BIMs) use an integral representation of the solution, namely the solution is defined by an integral of a suitable potential over the interface. 
If the domain boundaries are fixed and suitable parameterizations are available, BIMs can be a natural choice to solve exterior problems. 

In this paper, we present an implicit boundary integral algorithm for solving the following free-boundary problem
\begin{equation} \label{ms}
\begin{cases}
\Delta u(\mathbf{x},t) = 0 & \mathbf{x} \in \mathbb R^m \setminus \Gamma_{t},~m=2~\textrm{or}~3, \\
u(\mathbf{x}) = - \kappa(\mathbf{x}) & \mathbf{x} \in \Gamma_{t}, t \geq 0,  \\
v_{n} =  -  \left [ \frac{\partial u}{\partial \mathbf{n}} \right ]_{\Gamma_{t}} & \mathbf{x} \in \Gamma_{t}, t \geq 0, \\
\lim_{|\mathbf{x}| \rightarrow \infty} |u(\mathbf{x})|< \infty, & m=2, \\
\lim_{|\mathbf{x}| \rightarrow \infty} |u(\mathbf{x})| = 0, & m=3.\\
\Gamma_{0} = \partial \Omega_{0}^{-}. & \\
\end{cases}
\end{equation}
This problem is sometimes referred to as the Mullins-Sekerka problem \cite{mullins_sekerka63} which was first studied as a solidification process.  In general, this problem models the solidification or liquidation of materials of negligible specific heat, where the function $u$ represents the scaled temperature and $\Gamma_{t}$ is the boundary of the bounded domain $\Omega_{t}$, namely $\Gamma_{t} = \partial \Omega_{t}$ at time $t$. The boundary condition $u = - \kappa$ comes from the Gibbs-Thomson relation and is related to the surface tension effect. In two dimensions, the boundary condition at infinity (assuming that the temperature is uniformly bounded) produces an isolated boundary in which the flux $\pder{u}{n}$ is zero at infinity. In three dimensions, the boundary condition at infinity is an isothermal boundary condition which models the liquidation of material. 
  
In this setup, $\Omega_{t}$ is the solid region and $\Omega_{t}^{c}$ is the liquid region. $\left [ \frac{\partial u}{\partial \mathbf{n}} \right ]_{\Gamma_{t}}$ is the jump in the normal derivative of the function $u$ across the interface $\Gamma_{t}$ defined as the sum of the outward normal derivatives of $u$ from each side of $\Gamma_{t}$. 
More precisely,  
\begin{equation}
\left[\frac{\partial u}{\partial \mathbf{n}}\right]:= \left (\frac{\partial u}{\partial \mathbf{n^{+}}} \right )_{-} - \left (\frac{\partial u}{\partial \mathbf{n^{+}}} \right )_{+},
\end{equation}
where $+$ and $-$ refer to outside and inside respectively. 

 We observe that the two dimensional Mullins-Sekerka problem has nontrivial stationary solutions. For example, a collection of disjoint circles with the same radius $R$ will be stationary under the Mullins-Sekerka dynamics.  To see this, we note that in this case, both the solution inside and outside will be constant equal to $-\frac{1}{R}$, and thus $v_{n} = 0$.

This free boundary problem has interesting properties concerning volume and surface area. 
See for example  \cite{conti_niethammer_otto06}, where coarsening rates in Mullins-Sekerka evolution in bounded domains with periodic conditions are studied. 

In the following, we describe some notable differences in evolutions defined in two and three dimensions.
Let $A(t)$ denote the surface area of $\Gamma_{t}$ at time $t$ and $V(t)$ denote the volume of $\Omega_{t}$ at time $t$. Then we have 
\begin{align*}
\frac{dV(t)}{dt} & = - \int_{\Gamma_{t}} v_n ds \\
& =   \int_{\Gamma_{t}} \left [ \frac{\partial u}{\partial n} \right ]  ds \\
& =  \int_{\Omega_{t}} \Delta u dx -  \int_{\Gamma_{t}} \left  ( \frac{\partial u}{\partial \mathbf{n}^{+}}\right )_{+}ds.
\end{align*}
The first term is zero because $u$ is harmonic inside $\Omega_{t}$. In two dimensions, the second term turns out to be zero because of the compatibility condition  \eqref{comp_cond2D} (see Appendix \ref{comp_cond}), and thus the total volumeenclosed by $\Gamma_t$  is conserved. This compatibility condition is not valid in three dimensions and hence the motion is not volume preserving in three dimensions. 


For the area, we have the following:
\begin{align*}
\frac{dA(t)}{dt} & =   - \int_{\Gamma_{t}} \kappa v_{n} ds  \\
& = - \int_{\Gamma_{t}} u  \left [ \frac{\partial u}{\partial n} \right ] ds \\
& =  \int_{\Gamma_{t}} u \left (\nabla u \right )_{+}  \cdot  \mathbf{n}^{+}ds - \int_{\Omega_{t}} |\nabla u|^2 dx - \underbrace{\int_{\Omega_{t}} \Delta u dx }_{=0}.
\end{align*}
In two dimensions, we  can use equation~\eqref{eq_forlength} derived in Appendix~\ref{comp_cond} to
obtain
\begin{equation}
\frac{dA(t)}{dt} =  -\int_{\mathbb R^2 \setminus \Gamma_{t}}  |\nabla u|^2 dx,
\end{equation}
namely the evolution in two dimensions is area preserving while decreasing the length of the boundary.

Another interesting property of the Mullins-Sekerka dynamics is that the two dimensional motion induced by \eqref{ms} does not preserve convexity, see e.g. \cite{bates_chen_deng95}.

In this paper, we solve \eqref{ms} using an implicit boundary interface method, namely we define $\Gamma_{t}$ implicitly through a level set function $\varphi$, namely at any time $t \geq 0$, 
$$
\Gamma_{t} = \left \{ \mathbf{x} \in \mathbb R^m : \varphi(\mathbf{x},t) = 0\right \}.
$$
The algorithm operates by alternatively solving for the solution $u$ while maintaining the interface fixed, and then propagating the interface according to the normal velocity obtained from the jump in the normal derivative of $u$ across the interface. 
In the present work, the solution of the boundary value problem given by the first three equations in \eqref{ms} is solved using an implicit boundary integral formulation \cite{kublik_tanushev_tsai13, kublik_tsai16} that allows the boundary integral to be rewritten as an integral over $\mathbb R^m$.

\section{Integral equations for Laplace's equation} 
\label{int_eq_laplace}

We present below the boundary integral equation formulations most relevant to this paper and some of their useful properties. Throughout the paper, we consider $\Omega \subset \mathbb R^m$ to be a bounded set, $m=2$ or $3$ and the boundary  $\partial \Omega = \Gamma$ to be a disjoint collection of closed compact $C^2$ interfaces. We also denote by $\mathbf{n^{+}} $ and $\mathbf{n}^{-}$ the outward and the inward unit normal vectors respectively. The exterior is defined to be the unbounded domain. 
 
We consider the Dirichlet problem (interior and exterior) in an $L+1$ connected region (bounded or unbounded) with $L$ holes labeled $\left \{ S_{i} \right \}_{i=1}^{L}$ each having boundary $\Gamma_{i} = \partial S_{i}$:
\begin{equation}
\label{interior}
\begin{cases}
\Delta u (\mathbf{x}) = 0 & \mathbf{x} \in \Omega \\
u(\mathbf{x}) = f(\mathbf{x}), & \mathbf{x} \in \Gamma,
\end{cases}
\end{equation}
and 
\begin{equation}
\label{exterior}
\begin{cases}
\Delta u (\mathbf{x}) = 0 & \mathbf{x} \in \mathbb R^m \setminus \bar{\Omega}, \\
u(\mathbf{x}) = f(\mathbf{x}), & \mathbf{x} \in \Gamma, \\
\lim_{|\mathbf{x}| \rightarrow \infty} |u(\mathbf{x})|< \infty, & m=2 \\
\lim_{|\mathbf{x}| \rightarrow \infty} |u(\mathbf{x})| = 0, & m=3.
\end{cases}
\end{equation}

Since Dirichlet boundary conditions are imposed, we introduce an unknown density $\beta$ defined on $\Gamma$ and represent the solution of \eqref{interior} using the double layer potential formulation
$$
u(\mathbf{x}) = \int_{\Gamma} \beta(\mathbf{y}) \pder{\Phi(\mathbf{x},\mathbf{y})}{\mathbf{n_{y}^{+}}} dS(\mathbf{y}), ~~~ \mathbf{x} \in \Omega,
$$
where $\Phi$ is the fundamental solution of Laplace's equation defined as
\begin{equation} \label{fundsol}
\Phi(\mathbf{x},\mathbf{y}) = \left \{
\begin{array}{ll}
\frac{1}{2\pi} \ln |\mathbf{x}-\mathbf{y}| & \mbox{ for } m=2,\\
~\\
\ds -\frac{1}{m(m-2) \rho_{m}|\mathbf{x}-\mathbf{y}|^{m-2}} & \mbox{ for } m \geq 3,
\end{array}
\right.
\end{equation}
where $\rho_{m}$ is the volume of the unit ball in $\mathbb R^{m}$. 

In this context, the naive integral equation is ill-posed, both for the interior and the exterior problem, since there are $L$ nontrivial homogeneous solutions that span the nullspace, see e.g. \cite{folland76}. To alleviate this issue, we follow the approach first suggested by Mikhlin \cite{mikhlin57} and developed by Greenbaum etal \cite{greenbaum_greengard_macfadden93}. The procedure amounts to adding a linear combination  $ \sum_{i=1}^{L} A_{i} \Phi(\mathbf{x} - \mathbf{z_{i}})$ to the integral equation, where $\mathbf{z_{i}} \in S_{i}$, and the coefficients  $A_{i}$ play the role of Lagrange multipliers. 
These $L$ constants give $L$ degrees of freedom that are used to make sure the constructed solution satisfies the appropriate equation and boundary conditions. For the full derivation of the modified integral equation, we refer the reader to the work of Mikhlin and Greenbaum etal \cite{greenbaum_greengard_macfadden93, mikhlin57}. 

The interior problem is solved as follows:

\begin{enumerate}

\item Find the density $\beta$ and the constants $\left \{ A_{i} \right \}_{i=1}^{L}$ such that
\begin{equation}
\int_{\Gamma} \beta(\mathbf{y}) \pder{\Phi(\mathbf{x},\mathbf{y})}{\mathbf{n_{y}^{+}}} dS(\mathbf{y}) + \frac{1}{2} \beta(\mathbf{x}) +  \sum_{i=1}^{L} A_{i} \Phi(\mathbf{x} - \mathbf{z_{i}})  = f(\mathbf{x}), 
\mbox{ for } \mathbf{x} \in \Gamma. \label{inteq_interior1}
\end{equation}

\begin{equation}
\int_{\Gamma_{i}} \beta(\mathbf{y}) dS)\mathbf{y}) = 0,  \mbox{ for } 1 \leq i \leq L.
\label{inteq_interior2}
\end{equation}
\item Reconstruct the solution $u$ in $\Omega$ using the double layer potential formulation
$$
u(\mathbf{x}) = \int_{\Gamma} \beta(\mathbf{y}) \pder{\Phi(\mathbf{x},\mathbf{y})}{\mathbf{n_{y}^{+}}} dS(\mathbf{y}) +  \sum_{i=1}^{L} A_{i} \Phi(\mathbf{x} - \mathbf{z_{i}}) , ~~~~ \mbox{ for } \mathbf{x} \in \Omega.
$$
\end{enumerate}

Note that for simply connected $\Omega$ the constants $A_{i}$ are not needed since when $L=0$ the term $\sum_{i=1}^{0} A_{i} \Phi(\mathbf{x} - \mathbf{z_{i}}) = 0$ and equation \eqref{inteq_interior2} also removed.

For the exterior problem 
it is necessary to modify the integral equation in order to obtain an invertible system. It can be shown, e.g. \cite{kress99},  that if the kernel is modified using $\displaystyle \frac{1}{|\mathbf{x} - \mathbf{y}|^{m-2}}$, the system becomes well-posed. The exterior problem is therefore solved as follows:

\begin{enumerate}

\item Find the density $\beta$ and the constants $\left \{ A_{i} \right \}_{i=1}^{L}$ such that

\begin{equation}
\begin{aligned}
&\int_{\Gamma} \beta(\mathbf{y}) \left ( \pder{\Phi(\mathbf{x},\mathbf{y})}{\mathbf{n_{y}^{+}}} - \frac{1}{|\mathbf{x} - \mathbf{y}|^{m-2}} \right ) dS(\mathbf{y})- \frac{1}{2} \beta(\mathbf{x}) +  \sum_{i=1}^{L} A_{i} \Phi(\mathbf{x} - \mathbf{z_{i}})  = f(\mathbf{x}),\\
&\int_{\Gamma_{i}} \beta(\mathbf{y}) dS(\mathbf{y}) = 0,  \mbox{ for } 1 \leq i \leq L-1, ~~~\mbox{ for } \mathbf{x} \in \Gamma\\
&m=2:~~~\sum_{i=1}^{L} A_{i}  = 0, \\
&m=3:~~~ \int_{\Gamma_{L}} \beta(\mathbf{y}) dS(\mathbf{y}) = 0. 
\label{inteq_exterior}
\end{aligned}
\end{equation}

\item Reconstruct the solution $u$ in $\mathbb R^m \setminus \bar{\Omega}$ using the double layer potential formulation
$$
u(\mathbf{x}) = \int_{\Gamma} \beta(\mathbf{y}) \pder{\Phi(\mathbf{x},\mathbf{y})}{\mathbf{n_{y}}} dS(\mathbf{y}) +  \sum_{i=1}^{L} A_{i} \Phi(\mathbf{x} - \mathbf{z_{i}}) , ~~~~ \mbox{ for } \mathbf{x} \in R^m \setminus \bar{\Omega}.
$$
\end{enumerate}
The condition $\sum_{i=1}^{L} A_{i}  = 0$ necessary in the two dimensional case ensures that the solution is bounded. This condition is not needed in three dimensions.

\section{Overview of existing numerical methods}
In this section, we give a brief overview of numerical methods for solving integral equations. 
Boundary integral methods typically provide highly accurate numerical solutions. However, they require highly accurate parameterizations of the  boundaries, high order quadratures and specialized dense matrix solvers.

For each of the boundary integral equations obtained in the previous section, we need to solve a Fredholm equation of the second kind. In other words, we need to find a function $\gamma$ defined on $\Gamma$, such that
$$
q(x) = \int_{\Gamma} \gamma(y(s)) K(x,y(s)) ds + C_0 \gamma(x),
$$
where $C_0$ is a constant and $K$ is the normal derivative of the fundamental solution of Laplace's equation to $\Gamma$. To solve these equations
numerically it is necessary to discretize the above integrals. Three discretization methods are
typically used: the Nystr\"{o}m method \cite{atkinson97,nystrom30}, the collocation method \cite{atkinson97} and the
Galerkin method \cite{atkinson97,ciarlet78}. Each of these discretization methods leads to a discrete system of the
form
$$
(I + \mathtt{K}\Lambda) \gamma = q,
$$
where $I$ is the identity matrix, $\mathtt{K}$ is a dense matrix, $\Lambda$ is a
diagonal matrix (for example containing the quadrature weights of the
Nystr\"{o}m method), $\gamma$ is the vector of unknowns, and $q$ is a
known vector obtained from the boundary conditions. 

Since $\mathtt{K}$ is
dense this system is usually solved using an iterative procedure. In addition, low rank approximations may be constructed to improve the efficiency of the numerical solver. One very successful approach is the Fast Multipole
Method introduced by Greengard and Rokhlin in 1987
\cite{greengard_rokhlin87}. The use of hierarchical matrices \cite{born_grasedyck_hackbusch03}
to solve this dense system is also popular. 


In this paper, we use an exact integral formulation that allows the boundary integral to be rewritten exactly as a volume integral over a thin tubular neighborhood of the boundary. This enables us to perform the computations on a fixed gird regardless of the location of the interface. Should the interface evolve in time, all computations will be performed on the mesh that is used by the level set function at each time step. This makes the algorithm easy to implement for
evolving interfaces in two and three dimensions.

\section{Exact integral formulations using signed distance functions}

The exact integral formulation we use in this work was first proposed in \cite{kublik_tanushev_tsai13} and extended in \cite{kublik_tsai16}. This formulation allows the computation of integrals of the form 
\begin{equation}
\int_{\Gamma}v(\mathbf{x})dS,\label{intoverboundary}
\end{equation}
in the level set framework, namely when the domain $\Omega$ is represented
implicitly by a level set function.
%

In \cite{kublik_tanushev_tsai13}, with the choice of $\varphi=d_{\partial\Omega}$ 
being a signed distance function to $\Gamma$,  the integral \eqref{intoverboundary} is expressed as an average of integrals over nearby level sets of
$d_{\Gamma}$, where these nearby level sets continuously sweep a thin
tubular neighborhood around the boundary $\Gamma$ of radius
$\epsilon$. Consequently, \eqref{intoverboundary} is \emph{equivalent} to the volume integral shown on the right hand side below:
\begin{equation} 
\int_{\Gamma}v(\mathbf{x})dS=\int_{\mathbb{R}^{n}}v(\mathbf{x}^{*})J(\mathbf{x};d_{\Gamma})\delta_{\epsilon}(d_{\Gamma}(\mathbf{x}))d\mathbf{x},\label{intformulation}
\end{equation}
where $\delta_{\epsilon}$ is an averaging kernel, $\mathbf{x}^{*}$
is the
closest point  on $\Gamma$ to $\mathbf{x}$  and 
$J(\mathbf{x};d_{\Gamma})$  accounts 
for the change in curvature between the nearby level sets and the zero level set.  
We remark that the author in \cite{Miura-2015} used a similar formulation to study the heat equation defined in tubular neighborhoods of $\Gamma$.

We briefly give a justification for \eqref{intformulation}. Suppose that $\Gamma$ is a smooth closed hypersurface in $\mathbb{R}^m$ and assume that $\mathbf{x}$ is sufficiently close to $\Omega$ so that the closest point mapping
\[
\mathbf{x}^*=P_{\Gamma}(\mathbf{x})=\mathrm{argmin}_{y\in\Gamma} |\mathbf{x}-\mathbf{y}|
\]
 is continuously differentiable. 
 Then the restriction of $P_{\Gamma}$ to  $\Gamma_\eta$ is a diffeormorphism between $\Gamma_\eta$ and $\Gamma$, where $\Gamma_{\eta} \mathrel{\mathop:}=  \left \{ \mathbf{x} : d_{\Gamma} (\mathbf{x}) = \eta \right \}$.
As a result, it is possible to write integrals over $\Gamma$ using points on $\Gamma_\eta$ as:
\[
\int_{\Gamma} v(\mathbf{x}) dS = \int_{\Gamma_\eta} v(P_{\Gamma}(\mathbf{x}))J(\mathbf{x};\eta)dS,
\]
where the Jacobian $J(\mathbf{x},\eta)$ comes from the change of variable defined by $P_{\Gamma}$ restricted on $\Gamma_\eta$.
Averaging the above integrals respectively with a kernel, $\delta_\epsilon$, compactly supported in $[-\epsilon,\epsilon]$, we obtain
\[
\int_{\Gamma} v(\mathbf{x}) dS = \int_{-\epsilon}^{\epsilon} \delta_\epsilon(\eta) \int_{\Gamma_\eta} v(P_{\Gamma}(\mathbf{x}))J(\mathbf{x};\eta)dS~d\eta.
\]
Formula \eqref{intformulation} then follows from the coarea formula \cite{federer69} applied to the integral on the right hand side. Because the distance function is used, it is possible to obtain a closed form for the Jacobian $J(\mathbf{x},\eta)$. See \cite{kublik_tanushev_tsai13} and \cite{kublik_tsai16} for  proofs. 

%
We can now use this formulation to evaluate integral equations. Consider a general integral equation of the form
$$
\int_{\Gamma} \beta(\mathbf{y}) K(\mathbf{x},\mathbf{y}) dS(\mathbf{y}) + \lambda \beta(\mathbf{x}) = f(\mathbf{x}), ~~~~ \mbox{for } \mathbf{x} \in \Gamma,
$$
where $\lambda \in \mathbb R$ and $\Gamma$ is a closed $C^2$ hypersurface embedded in $\mathbb R^m$. We define for any function $u: \Gamma \mapsto \mathbb R$ its constant extension along the normals $\bar{u}$ as
$$
\forall \mathbf{x} \in \mathbb R^m, ~~~~~~~ \bar{u}(\mathbf{x}) = u(P_{\Gamma}(\mathbf{x})).
$$
Analytically, the formulation in \cite{kublik_tanushev_tsai13} gives us the following
$$
\int_{|d_{\Gamma}| \leq \epsilon} \bar{\beta}(\mathbf{y}) \tilde{K}(\mathbf{x},\mathbf{y}) \delta_{\epsilon}(d_{\Gamma}(\mathbf{y})) dy + \lambda \bar{\beta}(\mathbf{x}) = \bar{f}(\mathbf{x}),
$$
where $\tilde{K}(\mathbf{x},\mathbf{y}) = K(\mathbf{x},\mathbf{y}) J(\mathbf{x}; \Gamma)$ can be thought of as a weighted restriction.

On a uniform cartesian grid with grid spacing $h$, the above expression can be discretized as:  $\forall\mathbf{x}_{j} s.t. |d_{\Gamma}(\mathbf{x}_{j})| < \epsilon$

\begin{equation} \label{eq:KTT-linear-system}
{\sum_{k: |d_{\Gamma}(\mathbf{x}_{k})| < \epsilon} \bar{\beta}(\mathbf{x}_{k}) \tilde{K}(\mathbf{x}_{j}, \mathbf{x}_{k}) \delta_{\epsilon}(d_{\Gamma}(\mathbf{x}_{k})) h^m  + \lambda \bar{\beta}(\mathbf{x}_{j}) = \bar{f}(\mathbf{x}_{j}).}
\end{equation}
The number of equations in the corresponding system is the number of grid points trapped inside the tubular neighborhood $|d_{\Gamma}(\mathbf{x})| < \epsilon$.

\section{Algorithm for Mullins-Sekerka dynamics}

In this section, we highlight the main ideas of the algorithm and provide a discussion of the interesting and specific features of the algorithm.

\subsection{Algorithm}

There is variety of work that has been done on the computation of Mullins-Sekerka flows. Zhu etal \cite{zhu_chen_hou96} proposed a boundary integral method for computing the two space dimensional Mullins-Sekerka free boundary problem. In \cite{karali_kevrekidis09}, Karali etal used a set of ODE to simulate the interactions of circular bubbles under the Mullins-Sekerka flow both in two and three dimensions. In \cite{chen_merriman_osher_smereka97} the authors propose a level set based approach using finite difference to compute the two dimensional Mullins-Sekerka dynamics. Finally, in \cite{Barrett-Garcke-2010} the authors proposed a parametric finite element approach for two and three dimensions. In this paper, we propose an implicit boundary algorithm for simulating the resulting interfacial motion from Mullins-Sekerka dynamics on unbounded domains in both two and three dimensions. Our algorithm is based on the technique described in \cite{kublik_tanushev_tsai13}. However, because of the complexity of the Mullins-Sekerka dynamics, it is necessary to expand the basic method to more general configurations.

The overall structure of the proposed algorithm is defined by the general level set framework:
\begin{enumerate}
\item Given $d_{\Gamma_n}$ on the grid, solve the Laplace problems (5)-(6). Output: the density $\beta$ in $\{d_{\Gamma_n}(x_j)<\epsilon\}.$
\item Compute a suitable extension $\tilde{v}_n$ of $v_n$ in $\{d_{\Gamma_n}(x_j)<\epsilon\}.$
\item Solve the level set equation in $\{d_{\Gamma_n}(x_j)<\epsilon\}$  for $\Delta t$ amount of time 
\begin{equation}{\label{eq:levelset-PDE}
\varphi_t + \tilde{v}_n |\nabla\varphi|=0,
}
\end{equation}
using $\varphi(x,0)=d_{\Gamma_n}$ as initial condition.
\item Redistance: Given $\varphi(x,\Delta t)$, compute the signed distance function $d_{\Gamma_n}$ to the zero level set of $\varphi$.
\item Repeat the above steps until the desired final time. 
\end{enumerate}

In the following, we shall describe what we think are important details that are specific to the Mullins-Sekerka simulations, leaving the other steps to the standard level set literature.

\subsubsection{Dealing with multiply connected regions in Step 1}
\begin{figure}
\begin{centering}
\includegraphics{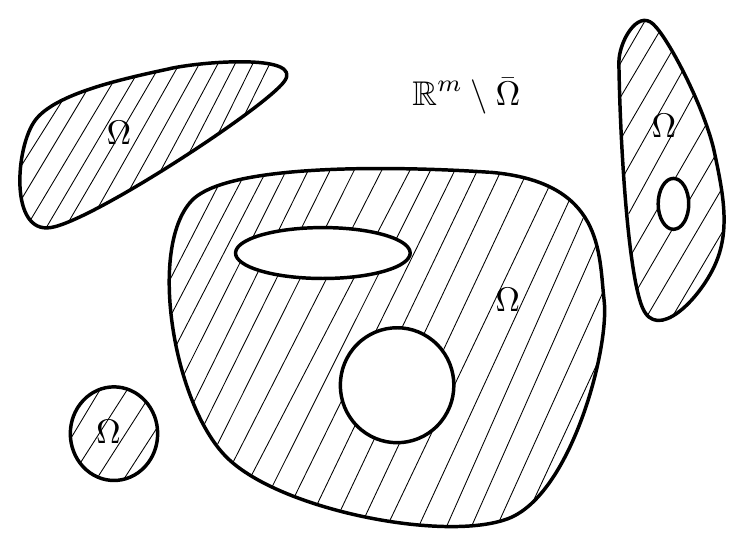}
\par\end{centering}
\caption{Non-simply connected domain.  
} 
\label{fig:components}
\end{figure}

One complexity is the fact that the regions can be multiply connected. 
Numerically, we first need to identify the connected components and how multiply connected they are before solving the correct system. See Figure~\ref{fig:components}  for an illustration of a typical configuration with simply and multiply connected components. To perform this task, we use a connected component labeling algorithm which was first proposed in the context of computer graphics. See Appendix~\ref{ccl}.

\subsubsection{Velocity extension in Step 2}

Another difficulty in the context of level set methods, is the velocity extension. In order to propagate an interface moving with normal velocity $v_{n}$, it is first necessary to extend $v_{n}(\mathbf{x},t)$ defined on $\Gamma_{t} \times [0,\infty)$  to $\tilde{v}_{n}(\mathbf{x},\tau)$ defined on $\mathbb R^m \times [0,\infty)$. This is typically done by solving the system
\begin{equation}
\label{velocity_extension}
\begin{cases}
\pder{\tilde{v}_{n}}{\tau} + sign(\varphi) \frac{\nabla \varphi}{|\varphi|} \cdot \nabla \tilde{v}_{n} = 0, \\
\tilde{v}_{n}(\mathbf{x},0) = v_{n}(\mathbf{x},t), \\
\end{cases}
\end{equation}
for a sufficient amount of time so that $\nabla \tilde{v}_{n} \cdot \nabla \varphi = 0$ in a neighborhood of $\Gamma_t$. This extension extends $v_{n}$ in a constant fashion along streamlines normal to the level sets of $\varphi$, which in practice is performed in a neighborhood of the interface, see \cite{cheng_tsai08} for a more extensive discussion on this topic.
However, this extension requires the knowledge of the normal velocity on the interface which in this case necessitates the evaluation of a hypersingular integral
$$
\pder{u}{\mathbf{n}_{\mathbf{x}}} (P_{\Gamma}(\mathbf{x})) = \int_{\Gamma} \frac{\partial ^2 \Phi}{\partial_{\mathbf{n_{\mathbf{x}}}} \partial_{\mathbf{n_{\mathbf{y}}}}} (P_{\Gamma}(\mathbf{x}), \mathbf{y}) \beta(\mathbf{y}) dS(\mathbf{y}).
$$

\begin{figure}
\begin{centering}
\includegraphics{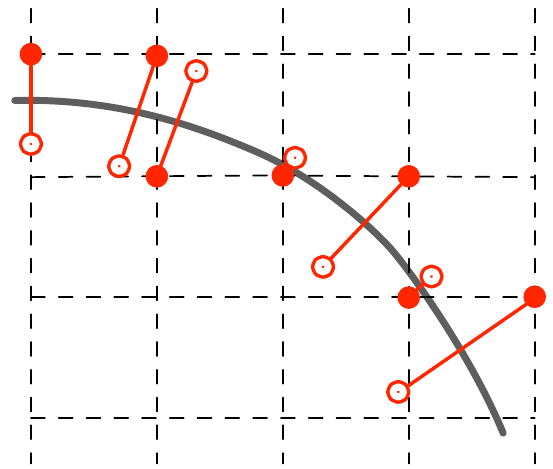}
\par\end{centering}
\caption{Velocity extension. The solid dots represent the grid points $\mathbf{x}$ and the empty ones their corresponding mirror points across the interface.  
} 
\label{fig:velocity-extension}
\end{figure}

We circumvent this issue by approximating the normal velocity $v_{n}$ at the grid points $\mathbf{x}$ closest to the interface $\Gamma_{t}$ (namely in a very narrow band around the interface taken to be a few grid points wide) using the equation
\begin{equation} 
\label{extension_mirror}
v_{n}(\mathbf{x}) \defeq sign(d_{\Gamma_{t}}(\mathbf{x})) \left (  \pder{u}{\mathbf{n}^{+}} (\mathbf{x}) -   \pder{u}{\mathbf{n^{+}}} \left (\mathbf{x} -  2 d_{\Gamma_{t}}(\mathbf{x}) \nabla d_{\Gamma_{t}}(\mathbf{x}) \right ) \right ), 
\end{equation}
where the normal derivatives $\pder{u}{\mathbf{n^{+}}}$ are computed using the interface integral
$$
\pder{u}{\mathbf{n}_{\mathbf{x}}} (\mathbf{x}) = \int_{\Gamma_{t}} \frac{\partial ^2 \Phi}{\partial_{\mathbf{n_{\mathbf{x}}}} \partial_{\mathbf{n_{\mathbf{y}}}}} (\mathbf{x}, \mathbf{y}) \beta(\mathbf{y}) dS(\mathbf{y}),
$$
where $\mathbf{x} \notin \Gamma_{t}$.  The idea behind equation~\eqref{extension_mirror} is to compute the jump in the normal derivative of $u$ at a point $\mathbf{x}$ close to $\Gamma_{t}$ by using it symmetrical point (or mirror point) $ \mathbf{x} -  2 d_{\Gamma_{t}}(\mathbf{x}) \nabla d_{\Gamma_{t}}(\mathbf{x})$ with respect to $\Gamma_{t}$. This is illustrated in Figure~\ref{fig:velocity-extension}. 

Note that the extension $\tilde{v}_{n}$ in \eqref{extension_mirror} is continuous and agrees with $v_{n}$ on $\Gamma_{t}$.  
We point out that this step can be interpreted as the initialization step for the PDE in equation~\eqref{velocity_extension}.  Thus, once the normal velocity $\tilde{v}_{n}$ is obtained in this very narrow band, we extend it to the next layer of grid points using the standard extension \eqref{velocity_extension}.

Finally the normal velocity may tend to $\infty$ as parts of $\Gamma_t$ vanish in finite time. Numerically, we regularize the computed normal velocity imposing a maximum normal velocity $v_{\infty}$. This means we replace $v_n$ by $\max(v_n,v_\infty)$.

\subsubsection{Literature on redistancing and involving the level set equations}

To maintain the completeness of this exposition, we briefly mention some established literature on the level set equation~\eqref{eq:levelset-PDE} and the construction of signed distance functions (Step 4 in the algorithm).

The level set equation \eqref{eq:levelset-PDE}, first proposed by Osher and Sethian in their seminal paper \cite{osher_sethian88}, convects the values of the level set function $\varphi$ with the velocity field $\tilde{v}_{n}$ which has been extended off of $\Gamma_{t}$. Typically, the extension is performed in a neighborhood of the interface since we are only interested in the location of the interface. Such techniques are referred to as local level set methods\cite{adalsteinsson_sethian95,peng_merriman_osher_zhao_kang99}.

The construction of signed distance functions (or distance reinitialization) is a procedure that replaces a general level set function by the signed distance function $d(\mathbf{x},t)$ which is the value of the distance from $\mathbf{x}$ to $\Gamma_{t}$ taken to be positive inside and negative outside (or vice versa). This assures that $\varphi$ does not become too flat nor too steep near the interface. To construct the signed distance function, one needs to find a function $d$ such that 
$$
|\nabla d| = 1 ~\mbox{almost everywhere},  ~~~~~~ \mbox{subject to } \left \{ \mathbf{x} : \varphi(\mathbf{x},t) = 0 \right \} = \left \{ \mathbf{x} : d(\mathbf{x},t) = 0 \right \}.
$$
Such construction can be performed very efficiently using fast sweeping or fast marching algorithms. See e.g. \cite{cheng_tsai08,russo_smereka00,sethian96,tsai_cheng_osher_zhao03,tsitsiklis95}.

\section{Numerical simulations}

In this section, we show the results of some numerical simulations using the algorithm. The aim of the following simulations is to show the computational properties of the algorithm.  
We discretize the integral equations \eqref{inteq_interior1}, \eqref{inteq_interior1} and \eqref{inteq_exterior} using finite differences and compute them using simple Riemann sums over uniform grids. After discretization, we obtain a linear system of the form $Ax = b$ where the vector $x$ consists of the unknown density $\beta$ and the constants $A_{i}$. See  \eqref{inteq_interior1}, \eqref{inteq_interior2} and \eqref{inteq_exterior}.

To solve equation \eqref{velocity_extension}, we use an upwind  scheme with a WENO-3 scheme for the spatial derivatives, and   third order TVD-RK3 method for time. We propagate enough to ensure the velocity information is transmitted throughout an $\epsilon$ wide band of the boundary $\Gamma_{t}$. For the level set equation, we use Godunov scheme with WENO-$3$ discretization for the eikonal term, and TVD-RK3 method for time.

\subsection{Two dimensions}

\begin{figure}[h]
\includegraphics[angle =0, width = 6cm]{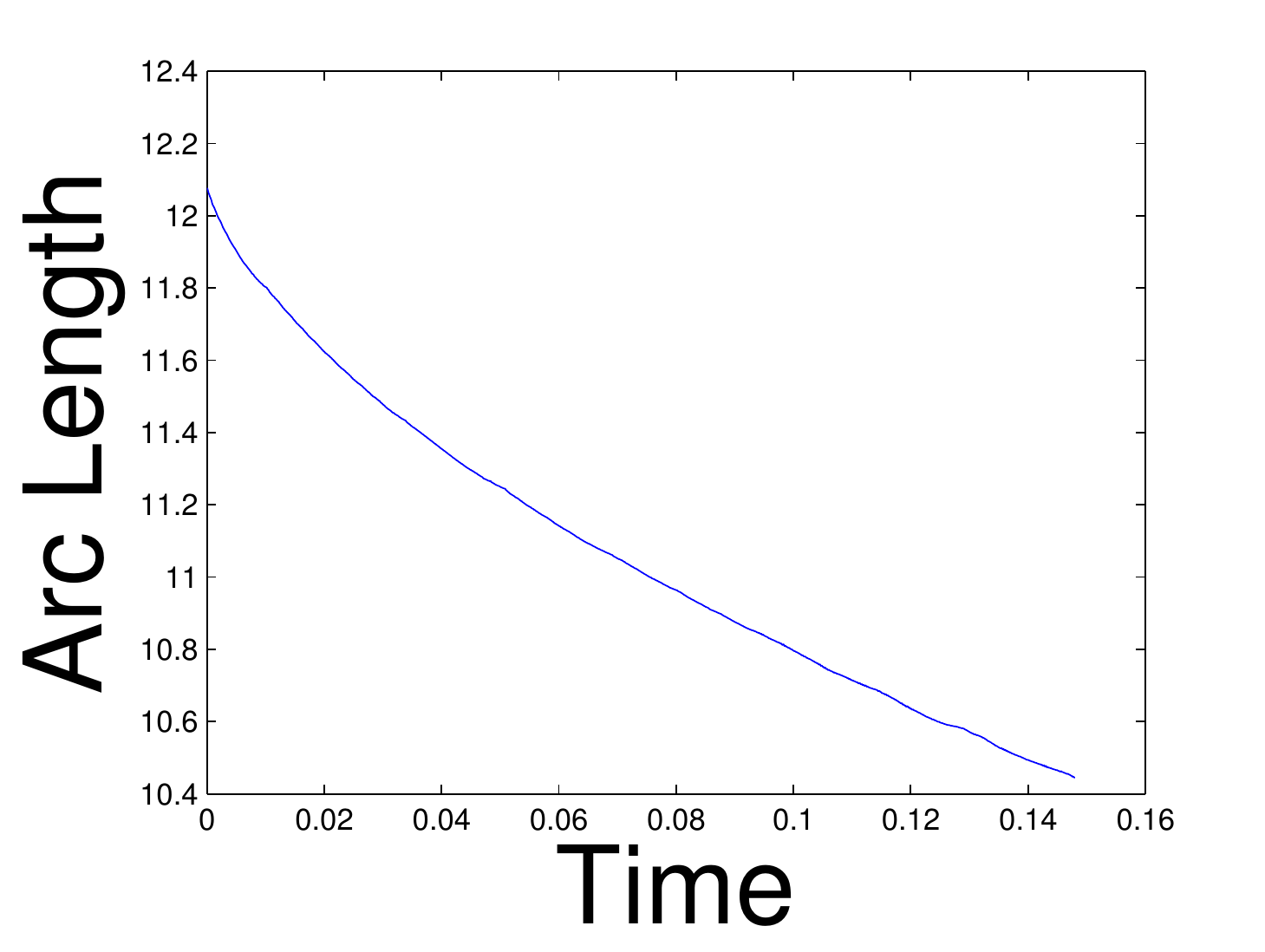}
\includegraphics[angle =0, width = 6cm]{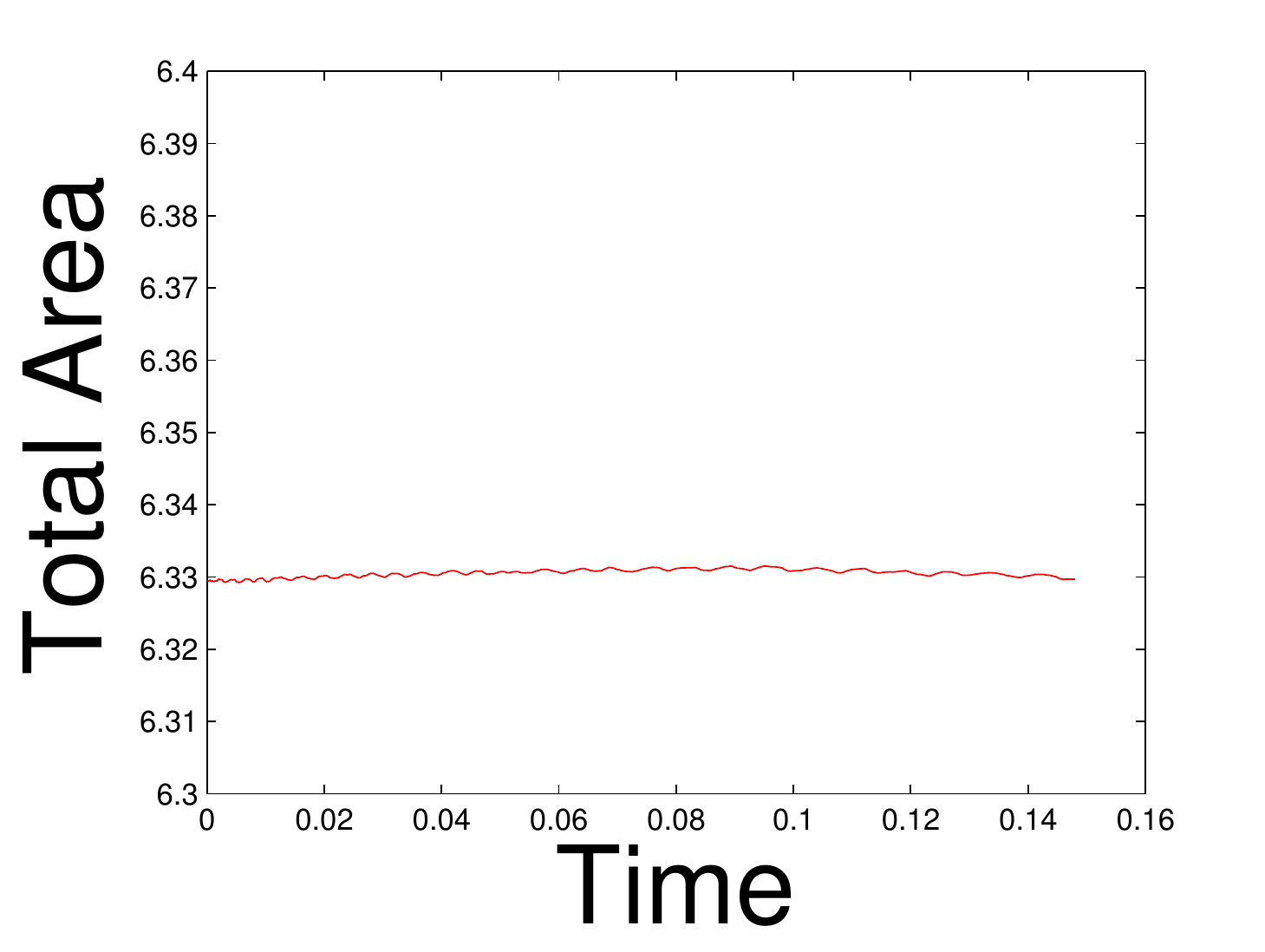}
\caption{Evolution of the perimeter and area of an ellipse under the Mullins-Sekerka model.}
\label{length-area-plot}
\end{figure}

We start by illustrating that the simulated two dimensional Mullins-Sekerka does indeed preserve area and decreases length. In this test case, we computed the evolution of the perimeter and area of an ellipse under the Mullins-Sekerka dynamics. The plots for the perimeter and area versus time are shown in Figure~\ref{length-area-plot}.

\begin{figure}[h]
\includegraphics[angle =0, width = 4.5cm]{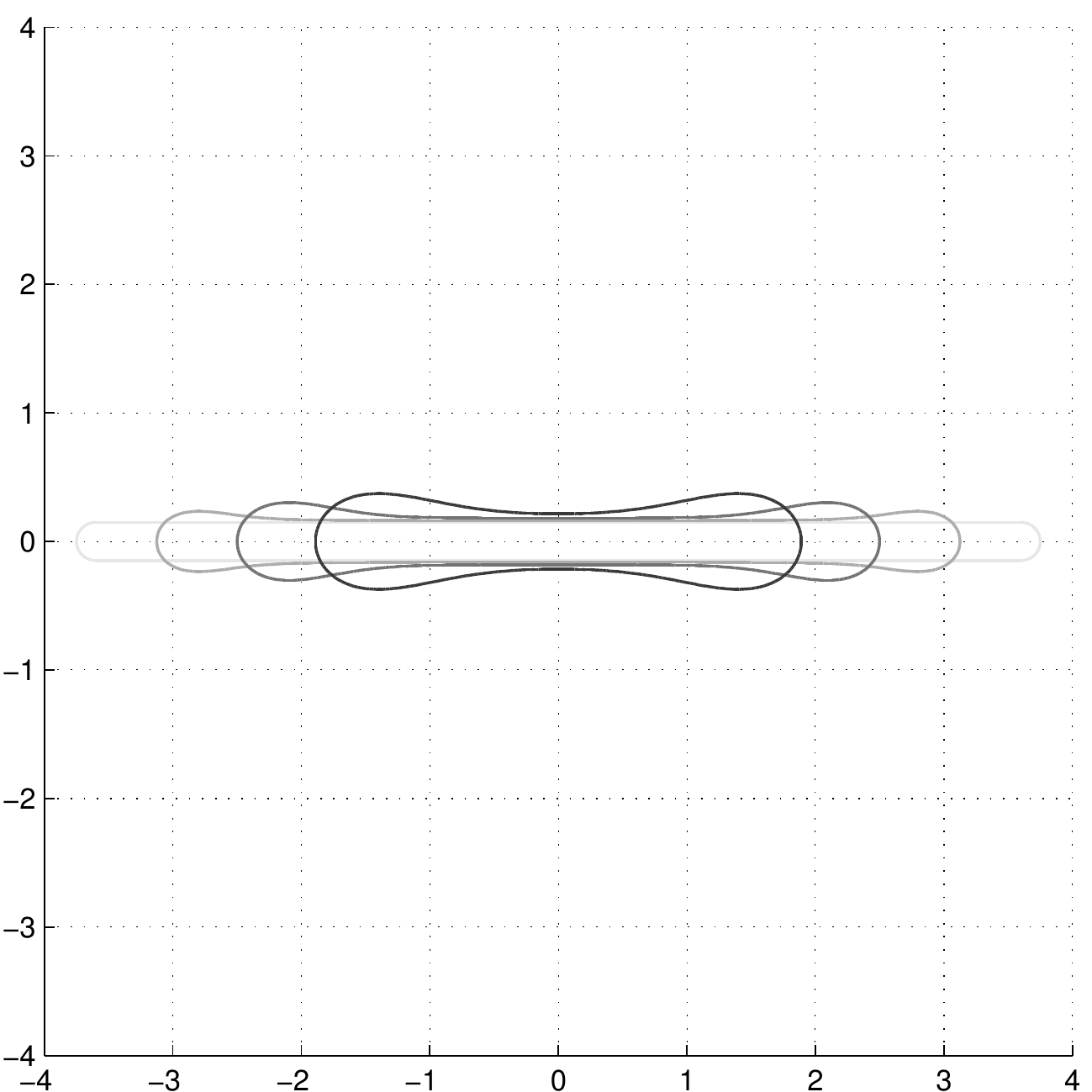}~~~~~~~~~~~~~~~~~~~~~~~~~~~~~~~
\includegraphics[angle =0, width = 4.5cm]{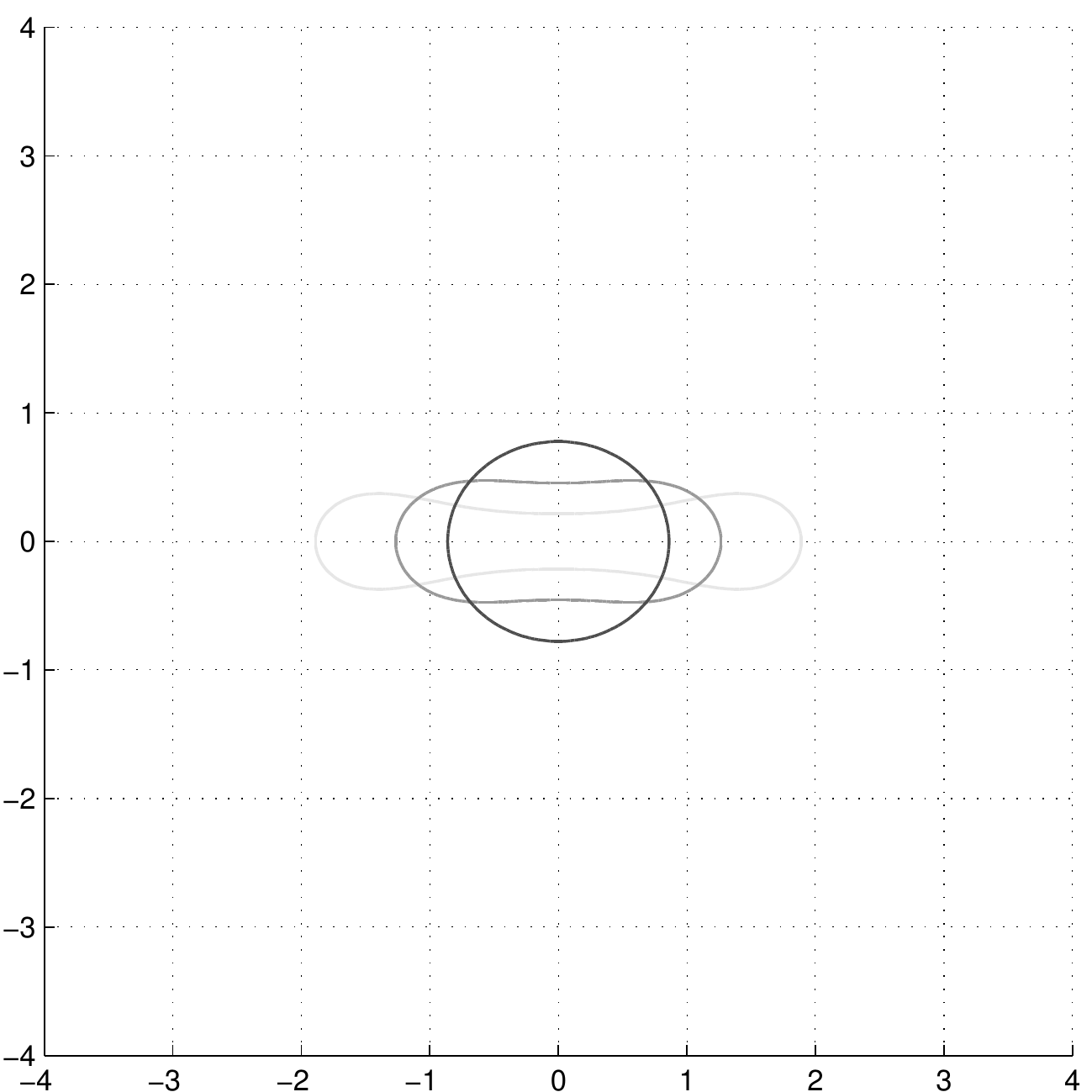}
\caption{Evolution of a thin tube under the two dimensional Mullins-Sekerka dynamics.}
\label{longtube-plot}
\end{figure}

As a second test problem, we simulated the evolution of an elongated tube to corroborate the fact that the two dimensional Mullins-Sekerka does not preserve convexity at all times during the evolution. This is illustrated in Figure~\ref{longtube-plot}. 

\begin{figure}[h]
\includegraphics[angle =0, width = 3.75cm]{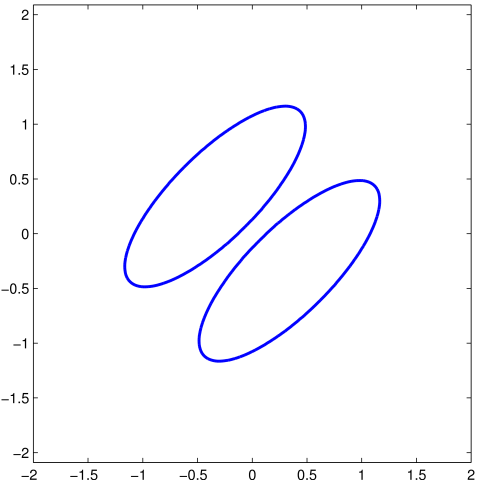}~~
\includegraphics[angle =0, width = 3.75cm]{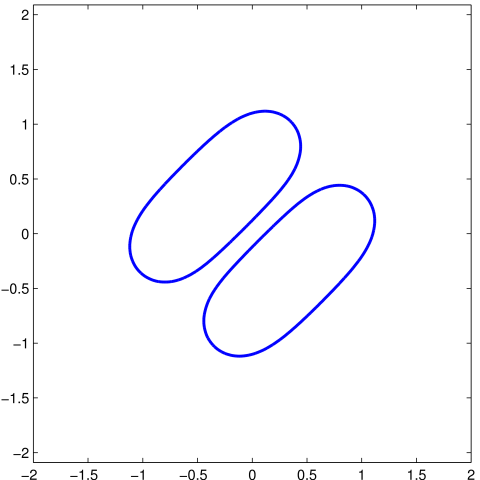}~~
\includegraphics[angle =0, width = 3.75cm]{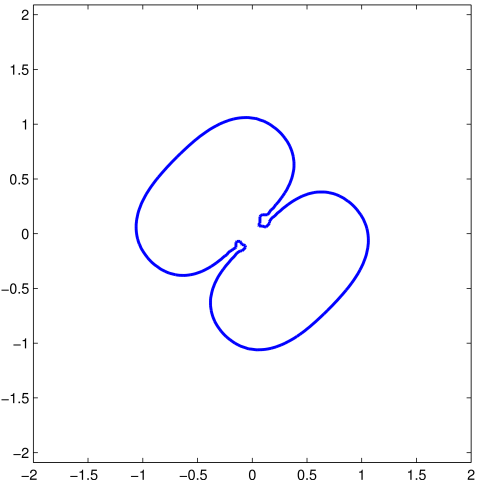} \\
\includegraphics[angle =0, width = 3.75cm]{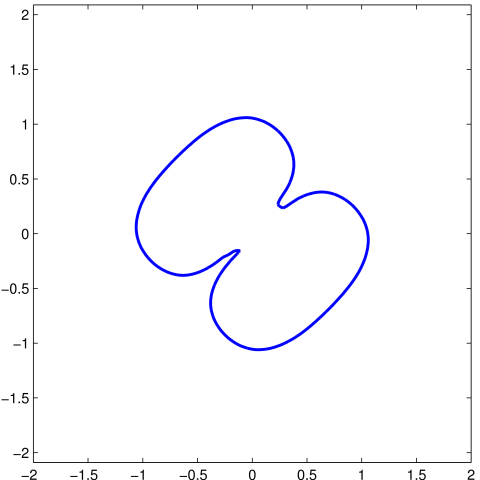}~~
\includegraphics[angle =0, width = 3.75cm]{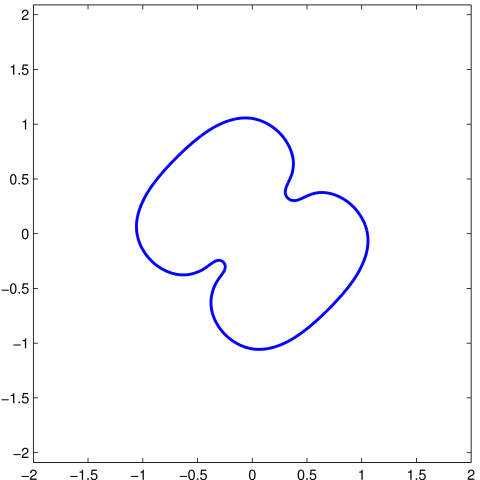}~~
\includegraphics[angle =0, width = 3.75cm]{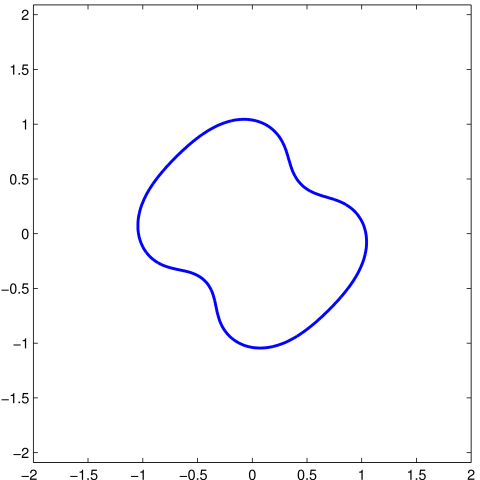} \\
\includegraphics[angle =0, width = 3.75cm]{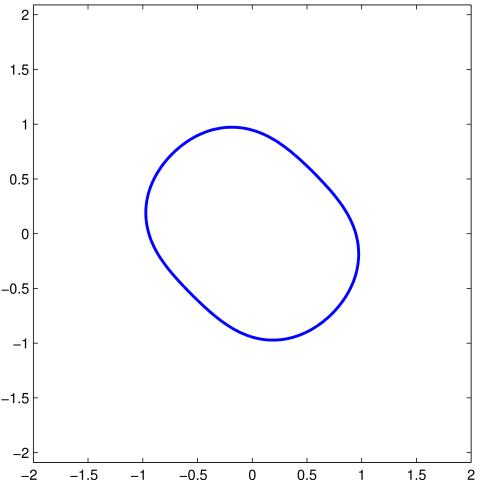}~~
\includegraphics[angle =0, width = 3.75cm]{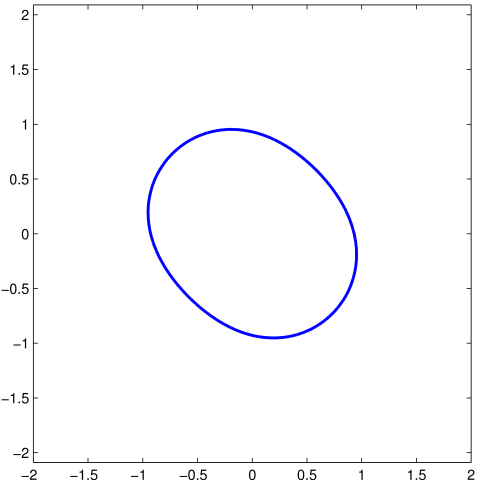}~~
\includegraphics[angle =0, width = 3.75cm]{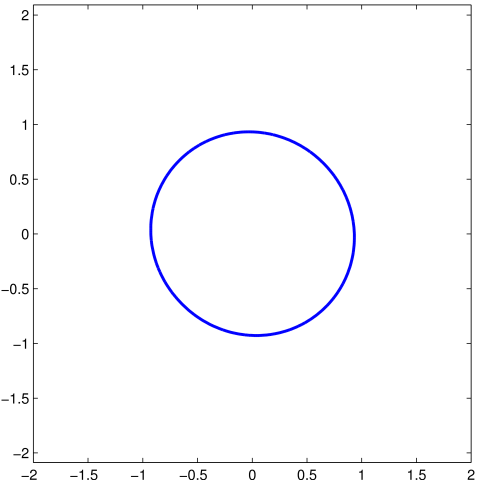}
\caption{Two ellipses merging in a two-dimensional Mullins-Sekerka flow simulation.}
\label{merging-plot}
\end{figure}

Finally, we demonstrate the capability of our algorithm to handle mergings of connected components. In this particular test example, we started with two ellipses that were close enough to finally merge into one connected component. Figure~\ref{merging-plot} shows the evolution of these two ellipses and Table~\ref{table-ellipse-merging} shows the numerical errors in the conservation of total area before and after merging.

\begin{table}[h] 
\begin{center}
\begin{tabular}{|c|c|c|c|}
\hline
h & $4/128$ & $4/256$ & $4/512$ \\
\hline
\hline
Initial Area & $2.6965$ & $2.69611$ & $2.695495$\\
\hline
Start Merge Area & $2.696655$ & $2.691792$ & $2.686751$ \\
\hline
End Merge Area & $2.804205$ & $2.724255$ & $2.694007$\\ 
\hline
Area Jump & $0.107550$ & $0.032463$ & $0.007256$ \\
\hline
Relative Area Error & $0.03988$ & $0.01206$ & $0.00270$\\
\hline
\end{tabular}
\end{center}
\caption{Merging time and area jump for the merging of two ellipses evolving under the two-dimensional Mullins-Sekerka dynamics.}
 \label{table-ellipse-merging}
\end{table}

\subsection{Three dimensions}
  
In these simulations, we assume that the far-field condition is $\lim_{|\mathbf{x}| \rightarrow \infty} |u(\mathbf{x})| = u_{\infty}$ which allows the solution to converge to a non-zero constant. This $u_{\infty}$ is referred to as the far-field environment which corresponds to the temperature at infinity. The dynamics of the Mullins-Sekerka model depend on this far-field environment and give rise to interesting behaviors that are illustrated in these simulations. While it is difficult to obtain rigorous analytical results on the general 3D Mullins-Sekerka model on unbounded domains with a general far-field environment $u_{\infty}$, we can simulate these different situations. We note that to obtain the so-called Mullins-Sekerka instability (or dendritic growth) it is necessary that $u_{\infty} < 0$. 

We observe that in section~\ref{int_eq_laplace}, we only provided a method for solving the three dimensional exterior problem with $u_{\infty}=0$. We thus briefly provide the equations for solving the \textbf{general 3D exterior problem}:

\begin{enumerate}

\item Find the density $\beta$ and the constants $\left \{ A_{i} \right \}_{i=1}^{L}$ such that
$$
\begin{cases}
\int_{\Gamma} \beta(\mathbf{y}) \left ( \pder{\Phi(\mathbf{x},\mathbf{y})}{\mathbf{n_{y}^{+}}} - \frac{1}{|\mathbf{x} - \mathbf{y}|^{m-2}} \right ) dS(\mathbf{y})- \frac{1}{2} \beta(\mathbf{x}) +  \sum_{i=1}^{L} A_{i} \Phi(\mathbf{x} - \mathbf{z_{i}})  = f(\mathbf{x}) - u_{\infty}, &
\\
\int_{\Gamma_{i}} \beta(\mathbf{y}) dS(\mathbf{y}) = 0,~~~1 \leq i \leq L,~~\mathbf{x} \in \Gamma. \\
\end{cases}
$$

\item Reconstruct the solution $u$ in $\mathbb R^m \setminus \bar{\Omega}$ using the double layer potential formulation
$$
u(\mathbf{x}) = \int_{\Gamma} \beta(\mathbf{y}) \pder{\Phi(\mathbf{x},\mathbf{y})}{\mathbf{n_{y}}} dS(\mathbf{y}) +  \sum_{i=1}^{L} A_{i} \Phi(\mathbf{x} - \mathbf{z_{i}}) + u_{\infty} , ~~~~ \mbox{ for } \mathbf{x} \in R^m \setminus \bar{\Omega}.
$$
\end{enumerate}

\begin{figure}[h]
\includegraphics[angle = 0, width = 3.85cm]{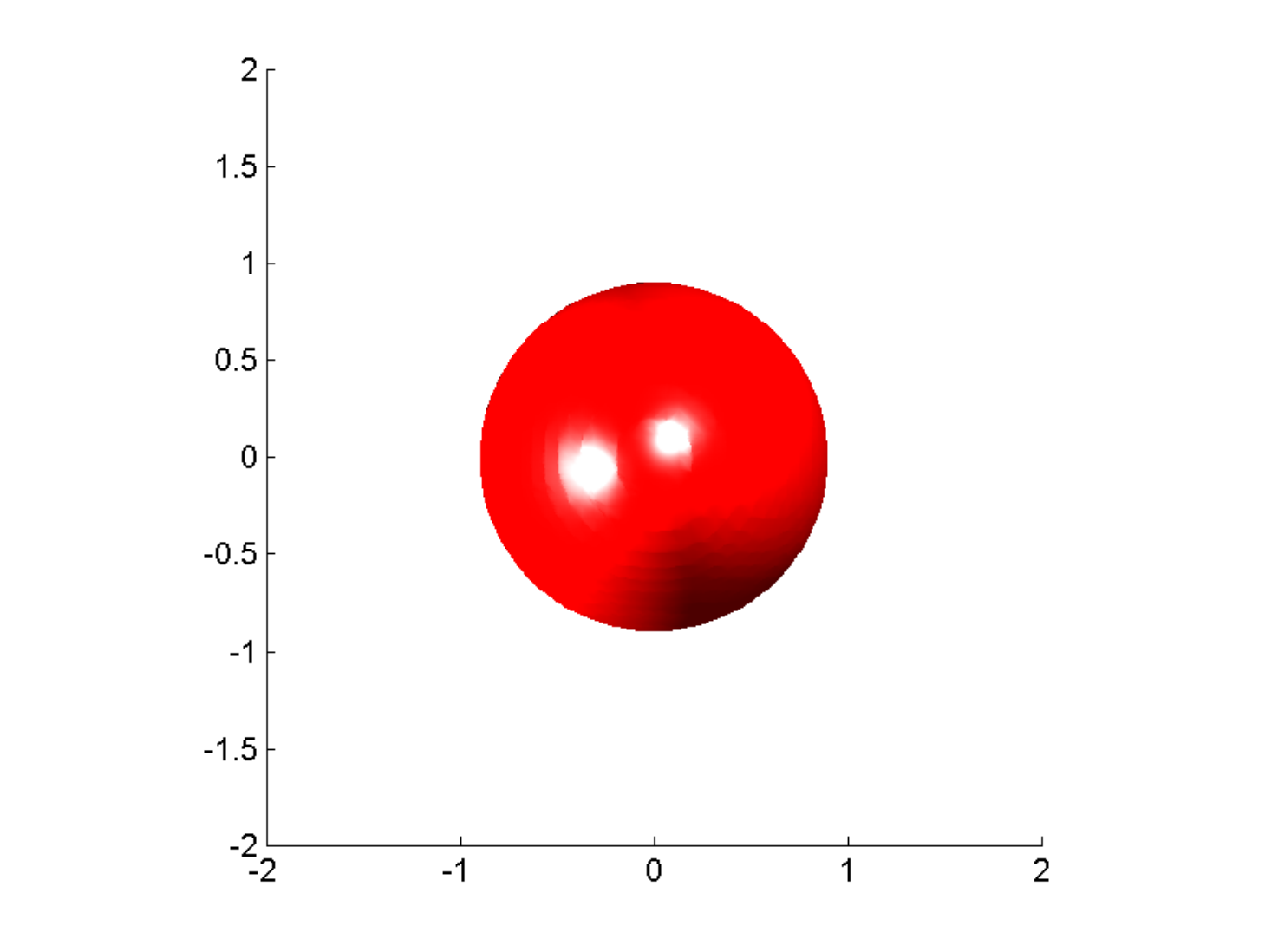}
\includegraphics[angle = 0, width = 3.85cm]{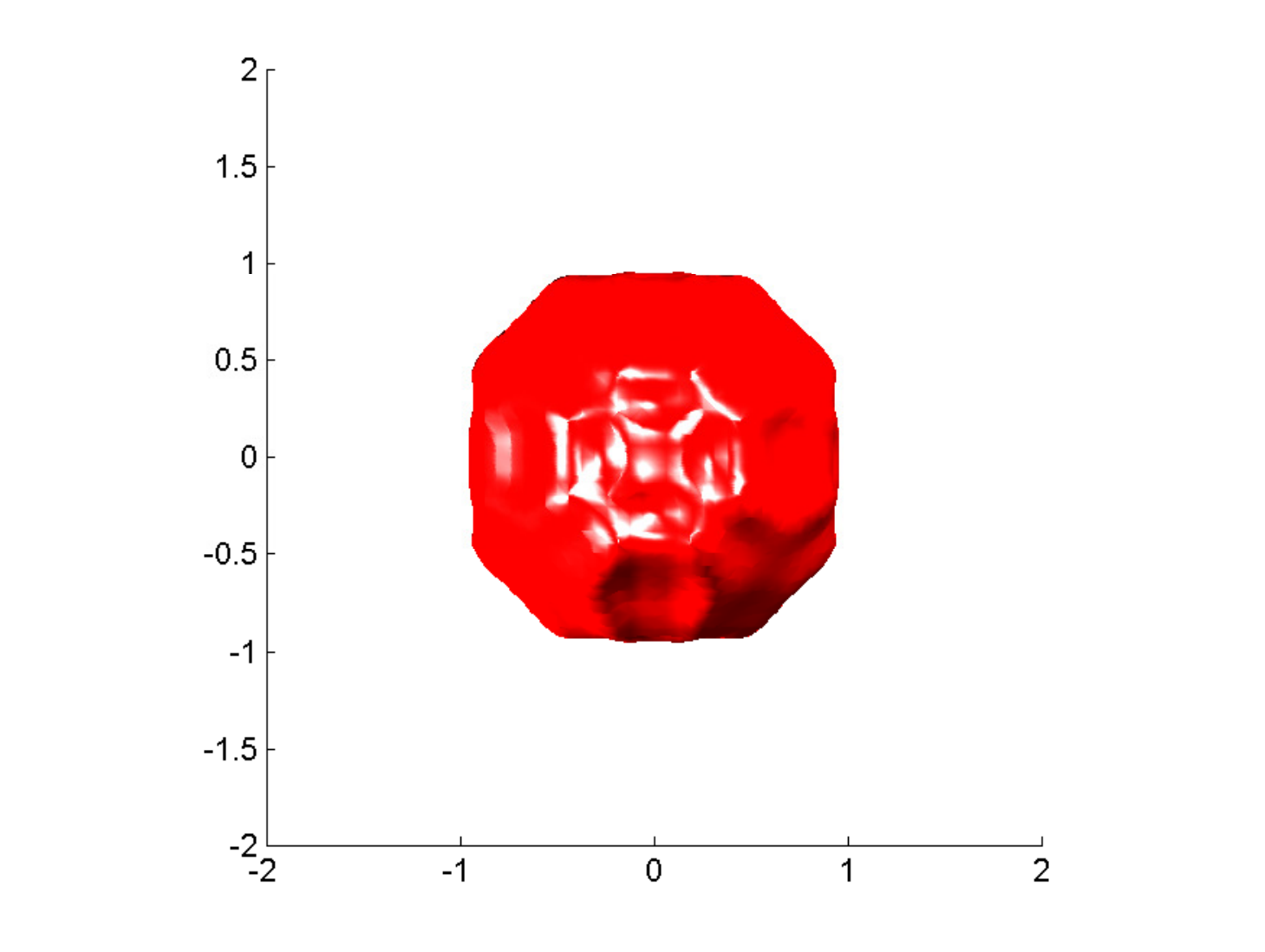}
\includegraphics[angle = 0, width = 3.85cm]{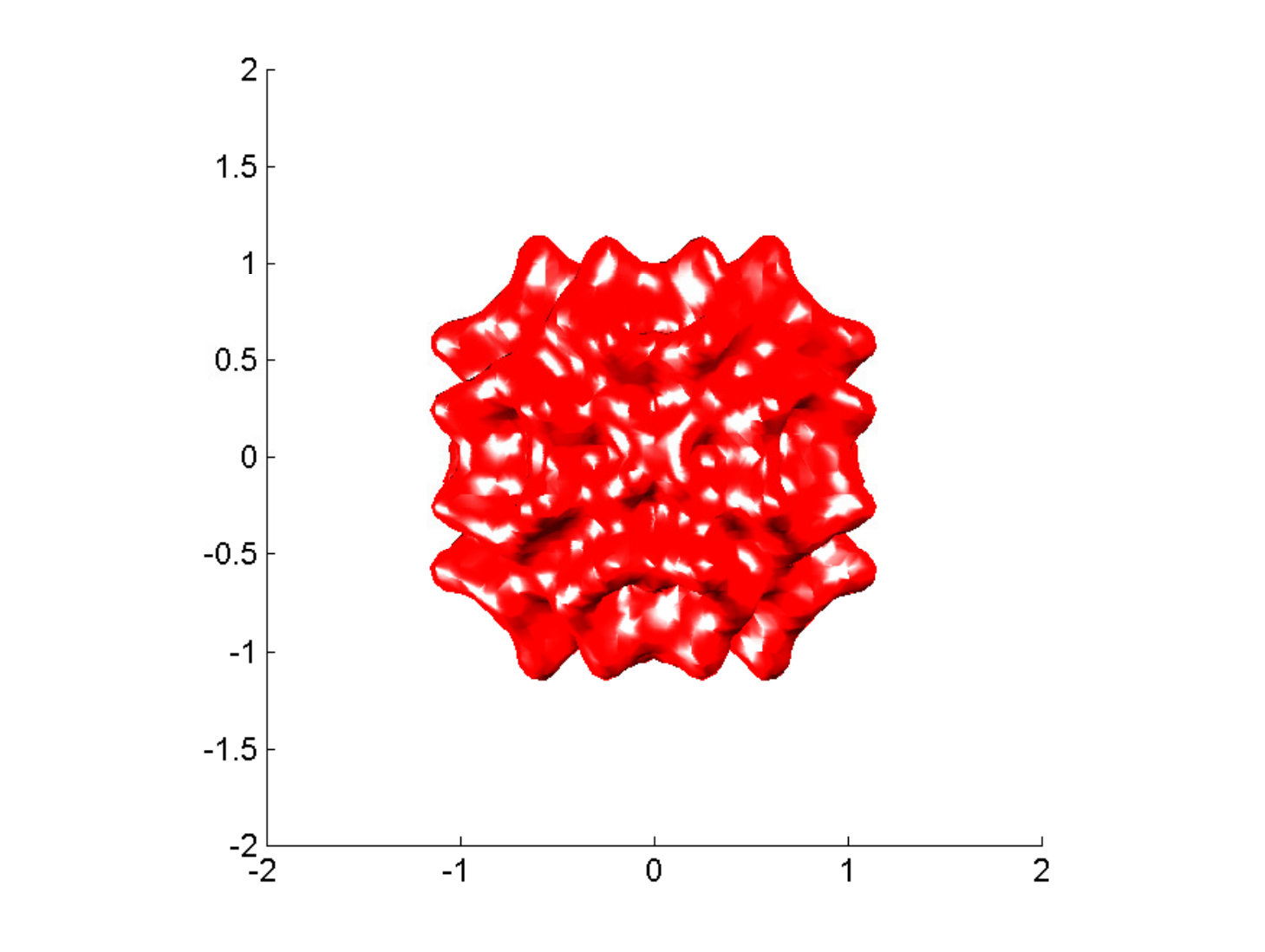}
\caption{Instability of dendritic growth in the 3D Mullins-Sekerka while maintaining crystal symmetry.}
\label{dendritic-plot}
\end{figure}

We focus on two features of the three dimensional Mullins-Sekerka dynamics. First, we look at the isothermal process and show that it depends on the far-field environment $u_{\infty}$. Second, we simulate a case where instabilities develop. See Figure~\ref{dendritic-plot}.

\begin{figure}[h]
\includegraphics[angle = 0, width = 3.85cm]{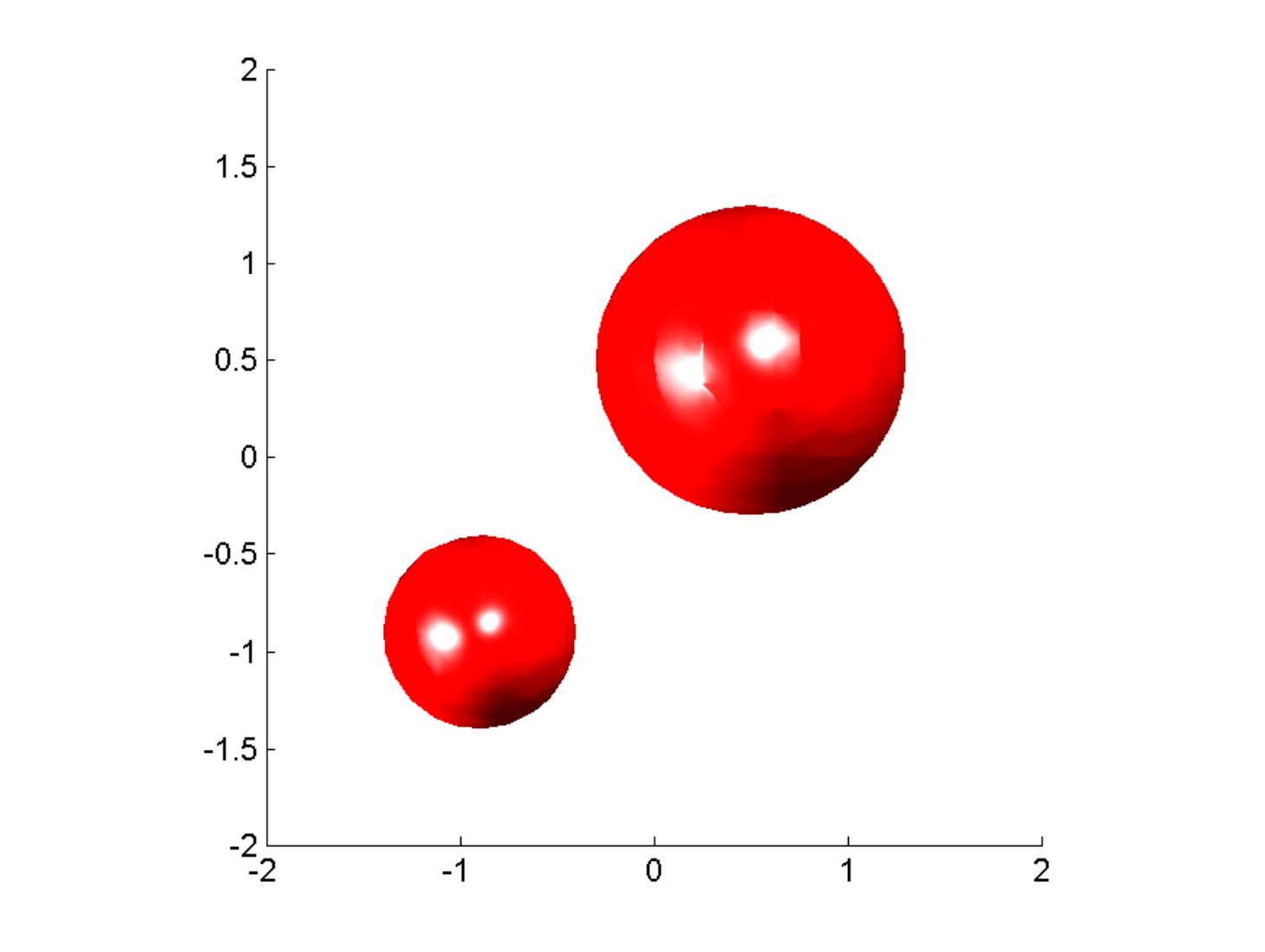}
\includegraphics[angle = 0, width = 3.85cm]{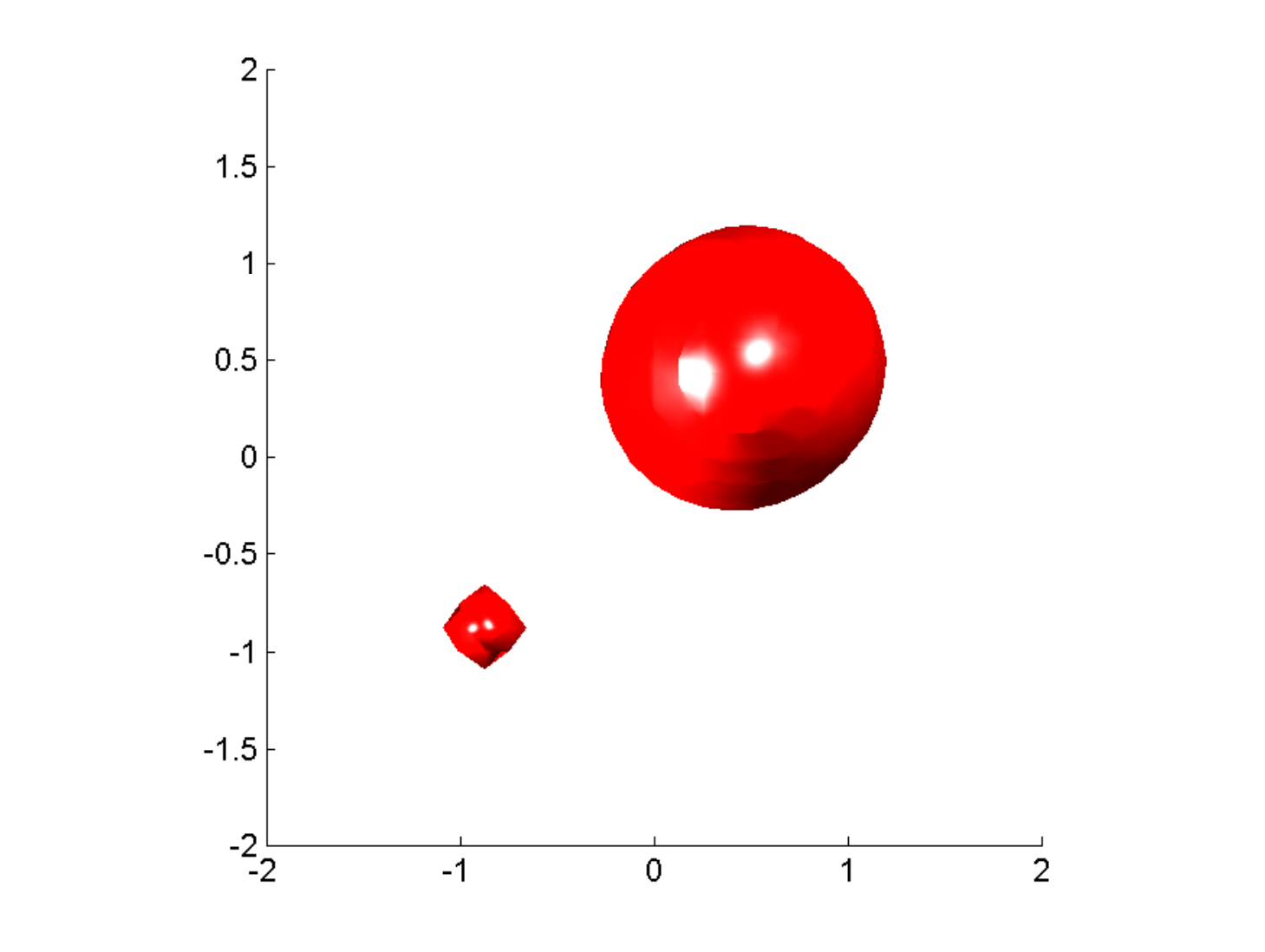}
\includegraphics[angle = 0, width = 3.85cm]{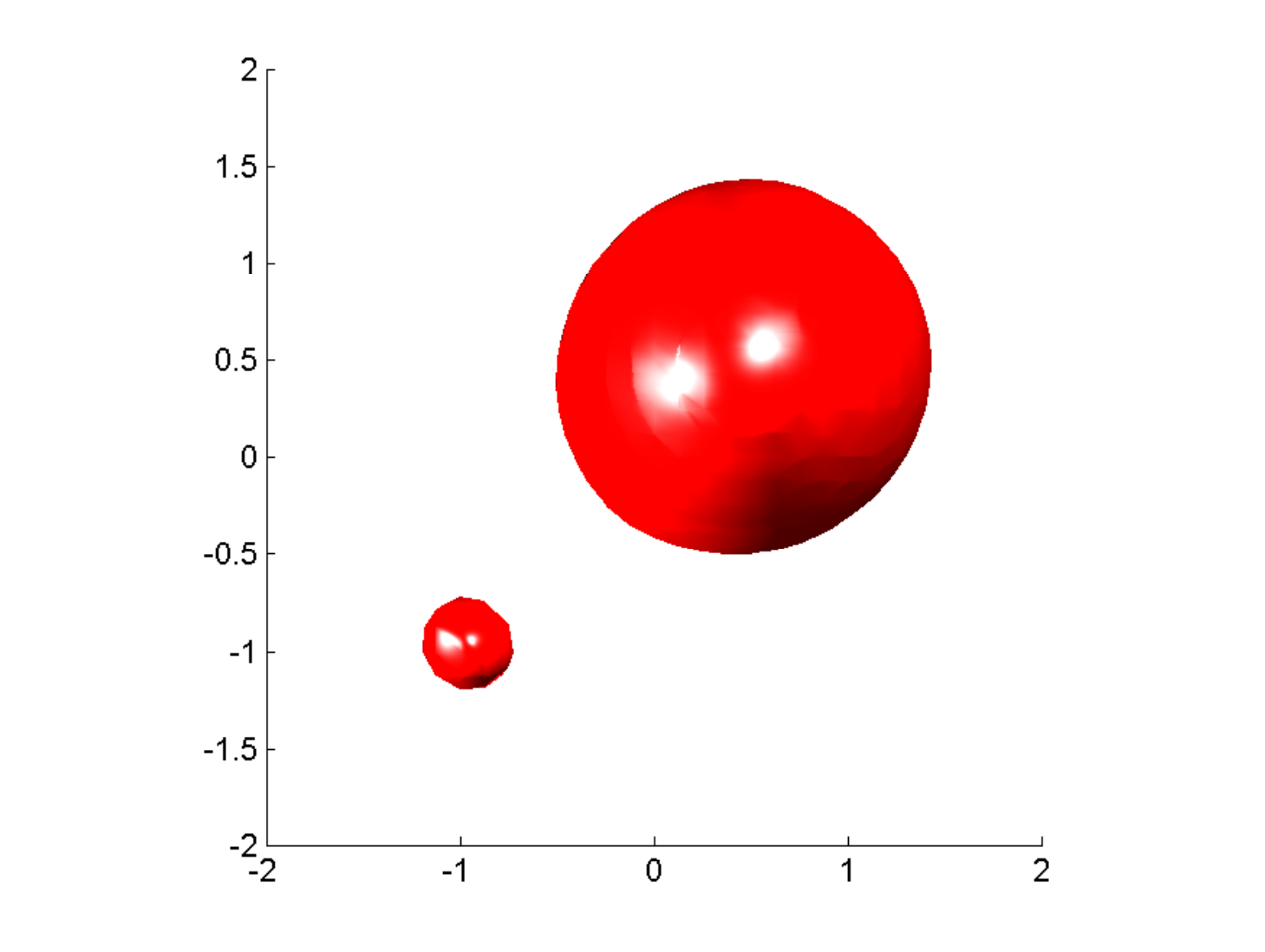}
\caption{Effect of the far-field value $u_{\infty}$ on the dynamics of two spheres. On the left plot, both spheres grow in a low temperature environment. In the middle, both shrink in a high temperature environment. On the right, the larger sphere grows and absorbes the mass of the smaller sphere.}
\label{effect-farfield-plot}
\end{figure}

In Figure~\ref{effect-farfield-plot} we simulate the dynamics of two initial spheres with different radii using various far-field value $u_{\infty}$.  If $u_{\infty}$ is small (i.e. low temperature), both spheres grow. If $u_{\infty}$ is large (high temperature), both spheres shrink (i.e. melt). If $u_{\infty}$ is set at phase transition point, we see that the larger sphere grows at the expense of the smaller one, as expected in the two dimensional case. 

\section{Conclusion}

We described an implicit boundary integral algorithm for computing the dynamics of the Mullins-Sekerka model on unbounded domains in two and three dimensions. We used a formulation described in \cite{kublik_tanushev_tsai13} and further investigated in \cite{kublik_tsai16} to compute integral equations using a signed distance function. We showed that our algorithm was able to handle topological changes, such as mergings of connected components, and was able to simulate the so-called Mullins-Sekerka instability (or dentritic growth). 

This paper demonstrates the versatility our computational method and in particular its capability to simulate on unbounded domains in three dimensions.  Coupled with fast multipole methods, this algorithm has the potential to be a powerful tool to better understand the dynamics of the Mullins-Sekerka model on  unbounded domains.
 
\begin{acknowledgement}
Kublik's research was partially funded by a University of Dayton Research Council Seed Grant and a Dr. Schraut Faculty Research Award. Chen's and Tsai's research is partially supported by Simons Foundation, NSF Grants
DMS-1318975, DMS-1217203, and ARO Grant No. W911NF-12-1-0519.
\end{acknowledgement}
\section*{Appendix}
\addcontentsline{toc}{section}{Appendix}

\subsection*{Connected component labeling} 
\label{ccl}

As explained in section~\ref{int_eq_laplace}, the system of equations to solve depends on many properties including the number of connected components of the region, whether the region is bounded or not (exterior vs interior), and the orientation of the region. 
In our numerical simulations, we adapted a technique called connected component labeling (CCL), see e.g.  
\cite{dillencourt_samet_tamminen92},
to find necessary topological information needed to obtain the correct set of equations. These include: 

\begin{enumerate}

\item The boundedness of the connected component $C_{i}$. This decides which formulation (interior vs exterior) to use.

\item The orientation of $C_{i}$ which determines the normal direction.

\item The total number of connected boundary components of the boundary of $C_{i}$, denoted by $\Gamma^{i} = \partial C_{i}$.

\item The selection of $\mathbf{z}_{i}$ in \eqref{inteq_interior1}, \eqref{inteq_interior2} and \eqref{inteq_exterior} for each component $\Gamma^j$. This involves 
	\begin{itemize}
	\item the separation of $\Gamma^{i}$ into boundary components $\Gamma^{i}_{j}$, each bounding a hole in the region except for the most exterior one. We denote the exterior boundary as $\Gamma^{i}_{0}$. 
	\item For each hole delimited by $\Gamma^{i}_{j}~,j \neq 0$, find a point $\mathbf{z}_{i}$ inside the hole that gives the least singular value of $\Phi(\mathbf{x}-\mathbf{z}_{i})$. This means that $\mathbf{z_{i}}$ should be as far from the boundary as possible.
	\end{itemize}
\end{enumerate}

Since we deal with closed interfaces, every grid point belongs to a unique connected component and has a well defined component label. Note that a connected component $C_{i}$ may belong to $\Omega$ or $\mathbf R^m \setminus \Omega$.  

We use the algorithm as follows:

\begin{itemize}

\item We label the unbounded component as $C_{0}$. This is the only true exterior component for the boundary integral formulations.

\item Each interior component ($d_{\Gamma}>0$, i.e. the solid phase) will have positive label $i$ and each exterior component ($d_{\Gamma}<0$ i.e. the liquid phase) will take negative label $-i$ (except for the unbounded one). 

\item Each boundary piece of each component will have label $j$. 

\item The point $\mathbf{z_{i}}$ in $C_{i}$ is chosen to be the point with largest signed distance function in absolute value.

\end{itemize}

We adopt the two-pass CCL algorithm with $2m$-connectivity in $\mathbb R^m$. 
This algorithm uses equivalence classes for labels: after the first pass, points in the same connected component may not have the same label, but the labels of points in the same connected component will be assigned to the same equivalent class in the second pass.
The root of the equivalence class denotes the smallest label (in absolute value) in the equivalence class. The scanning process works as follows:

\begin{enumerate}

\item First Pass

\begin{enumerate}

\item With begin with label $0$ which denotes the unbounded component $C_{0}$. 

\item At the current point scanned, we look at the sign of the signed distance function of its $2m$ neighbors already visited (in 2D, west and north of the current point)

\begin{itemize}

\item If none of the neighbors have the same sign as the current point, we create a new label (positive or negative, based on the sign of the distance function at the current point) and a new equivalence class. 

\item If only one neighbor has the same sign, we pick the label for the current point to be the root of that neighboring point's equivalence class. 

\item If more than one neighbor has the same sign as the current point, we pick the label for the current point to be the smallest root of the neighbors that have the same sign's equivalence classes. Furthermore, we combine the equivalence classes of the neighboring points that have the same sign since they are connected through the current point.

\end{itemize}

\item The largest root of all equivalence classes with a given sign (interior or exterior) will give the total number of connected components that have that sign. Thus, the sum is the total number of connected components.

\end{enumerate}

\item Second Pass

\begin{enumerate}

\item At each point, we assign its label to be the root of its equivalence class.

\item We update and store the points with largest absolute distance within each equivalence class. These are the points $\mathbf{z_{i}}$.

\item For each point $\mathbf{x}$ within the $\epsilon$ tubular neighborhood of the boundary $(|d_{\Gamma}|< \epsilon)$, we look for the root of its equivalence class (say $i$) and the root of its projection point's ($P_{\Gamma}(\mathbf{x})$) equivalence class (say $j$) that takes the opposite sign.  To obtain $j$, we look at the vertices of the cell $P_{\Gamma}(\mathbf{x})$ falls in and scan for the label with opposite sign. This step identifies the boundary piece $\Gamma^{i}_{j}$ and collects points within the $\epsilon$ neighborhood of the boundary $\Gamma^{i}_{j}$. 

\item We 
 store the total number of connected boundary pieces. 

\end{enumerate}

\end{enumerate}

\subsection*{Compatibility conditions}
\label{comp_cond}

Consider $\Omega$ an open bounded domain in $\mathbb R^m$ with smooth boundary $\partial \Omega = \Gamma$ for $m=2,3$. Let $R>0$ be sufficiently large so that $\bar{\Omega} \subset B_{R}(x)$ where $B_{R}(x)$ is the ball centered at $x \in \Omega$ with radius $R$.  Define the interior problem 
$$
\begin{cases}
\Delta u = 0 & \mbox{ in } \Omega \\
u = f & \mbox{ on } \Gamma,
\end{cases}
$$
and the exterior problem
$$
\begin{cases}
\Delta v = 0 & \mbox{ in } B_{R}(x) \setminus \bar{\Omega} \\
v = u & \mbox{ on } \Gamma \\
\mbox{ "far field condition"}.
\end{cases}
$$
We discuss what "far field condition" should be in the following.
Using Green's identity in $B_{R}(x) \setminus \bar{\Omega}$ for $x \in \Omega$, we obtain
\[
\begin{aligned}
0= \int_{B_{R}(x) \setminus \bar{\Omega}} &\left ( v(y) \Delta_{y} \Phi(x,y) - \Phi(x,y) \Delta_{y} v(y) \right ) dy \\
= &  \int_{\Gamma} \left (v(y) \pder{\Phi(x,y)}{\mathbf{n}_{y}}  - \Phi(x,y) \pder{v(y)}{\mathbf{n}_{y}}\right ) dS(y) \\
& - \int_{\partial B_{R}(x)} \left (v(y) \pder{\Phi(x,y)}{\mathbf{n}_{y}}  - \Phi(x,y) \pder{v(y)}{\mathbf{n}_{y}} \right ) dS(y) .
\end{aligned}
\]
Thus
\begin{equation}
\int_{\Gamma} \left (v(y) \pder{\Phi(x,y)}{\mathbf{n}_{y}}  - \Phi(x,y) \pder{v(y)}{\mathbf{n}_{y}} \right ) dS(y)  
=  \pder{\Phi(R)}{r} \int_{\partial B_{R}(x)} v ds - \Phi(R) \int_{\partial B_{R}(x)} \pder{v}{r} ds
\label{LHSRHS}
\end{equation}
with $\Phi(x,y) = \Phi(|x-y|) = \Phi(r)$. 
In addition, notice that since $v$ is harmonic in $B_{R}(x) \setminus \bar{\Omega}$, 
$$
\int_{\Gamma} \left ( \pder{v}{n} \right)_{+} ds = \int_{\partial B_{R}(x)} \pder{v}{r}ds.
$$

Now, since formula \eqref{LHSRHS} holds for any $R>0$,  we see that necessarily for $m=2$
\begin{equation}
\int_{\Gamma} \left ( \pder{v}{n} \right)_{+} ds=0,
\label{comp_cond2D}
\end{equation}
otherwise the second integral on the right hand side will not be bounded uniformly in $R$.
For $m=3$, it suffices that  the integral $\int_{\Gamma} \left ( \pder{v}{n} \right)_{+} ds$ be bounded uniformly in $R$, since $\Phi(R)\sim O(1/R)$.

Consequently 
"the far field condition" should be 
\[m=2:~~~\lim_{R\rightarrow\infty}|v|<\infty, \]
since $\pder{\Phi(r)}{r} |_{r = R} $ is $O(R^{-1})$ and $ds \sim  O(R)$, and
\[
m=3:~~~\lim_{R\rightarrow0}v=0,
\]
 since $\Phi(R)\sim{O}(R^{-1})$ and $ds\sim O(R^2)$.


Now applying the divergence theorem in the ``annulus" we obtain
\begin{align}
\int_{\partial B_{R}(x)} v \nabla v \cdot \mathbf{n}^{+} ds - \int_{\Gamma} v \nabla v \cdot \mathbf{n}^{+} ds & = \int_{B_{R}(x) \setminus \bar{\Omega}} \nabla \cdot \left ( v \nabla v \right ) dx \\
& = \int_{B_{R}(x) \setminus \bar{\Omega}} |\nabla v|^2 dx + \int_{B_{R}(x) \setminus \bar{\Omega}} \Delta v dx\\
& =  \int_{B_{R}(x) \setminus \bar{\Omega}} |\nabla v|^2 dx \label{eq_div_annulus}
\end{align}
since $v$ is harmonic in $\mathbb R^m \setminus \bar{\Omega}$. Now in two dimensions and by the previous argument, we know that $\lim_{R \rightarrow \infty} \pder{v}{n}|_{\partial B_{R}(x)} = 0$ which implies that for $m=2$
$$
\lim_{R \rightarrow \infty} \int_{\partial B_{R}(x)} v \nabla v \cdot \mathbf{n}^{+} ds = \lim_{R \rightarrow \infty} \int_{\partial B_{R}(x)} v \pder{v}{\mathbf{n}^{+}}ds = 0.
$$
Thus, taking $R$ to infinity in equation \eqref{eq_div_annulus}, we obtain
\begin{equation}
\int_{\Gamma} v \nabla v \cdot \mathbf{n}^{+} ds = - \int_{\mathbb R^2 \setminus \bar{\Omega}} |\nabla v|^2 dx.
\label{eq_forlength}
\end{equation}



%
%
%

\bibliographystyle{abbrv}
\bibliography{references}

\begin{thebibliography}{10}

\bibitem{adalsteinsson_sethian95}
D.~Adalsteinsson and J.~A. Sethian.
\newblock A fast level set method for propagating interfaces.
\newblock {\em J. of Comput. Phys.}, 118(2):269--277, 1995.

\bibitem{atkinson97}
K.~E. Atkinson.
\newblock {\em The Numerical Solution of Integral Equations of the Second
  Kind}.
\newblock Cambridge University Press, 1997.

\bibitem{atkinson_chandler90}
K.~E. Atkinson and G.~Chandler.
\newblock Boundary integral equation methods for solving laplace's equation
  with nonlinear boundary conditions: the smooth boundary case.
\newblock {\em Mathematics of Computation}, 55(191):451--472, 1990.

\bibitem{babuska70}
I.~Babu\u{s}ka.
\newblock The finite element method for elliptic equations with discontinuous
  coefficients.
\newblock {\em Computing}, 5:207--213, 1970.

\bibitem{Barrett-Garcke-2010}
J.~W. Barrett, H.~Garcke, and R.~Nürnberg.
\newblock On stable parametric finite element methods for the stefan problem
  and the mullins–sekerka problem with applications to dendritic growth.
\newblock {\em Journal of Computational Physics}, 229(18):6270 -- 6299, 2010.

\bibitem{bates_chen_deng95}
P.~W. Bates, X.~Chen, and X.~Deng.
\newblock A numerical scheme for the two phase {M}ullins-{S}ekerka problem.
\newblock {\em Electronic Journal of Differential Equations}, 1995.

\bibitem{bedrossian_teran10}
J.~Bedrossian, J.~J. von Brecht, S.~Zhu, E.~Sifakis, and J.~Teran.
\newblock A second order virtual node method for elliptic problems with
  interfaces and irregular domains.
\newblock {\em J. Comput. Phys.}, 229:6405--6426, 2010.

\bibitem{born_grasedyck_hackbusch03}
S.~B{\"{o}}rn, L.~Grasedyck, and W.~Hackbusch.
\newblock Hierarchical matrices.
\newblock Technical report, Max-Planck Institut fur Mathematik in den
  Naturwissenschaften, Leipzig, 2003.

\bibitem{chen_merriman_osher_smereka97}
S.~Chen, B.~Merriman, S.~Osher, and P.~Smereka.
\newblock A simple level set method for solving {S}tefan problem.
\newblock {\em J. Comput. Phys.}, 135, 1997.

\bibitem{cheng_tsai08}
L.-T. Cheng and Y.-H. Tsai.
\newblock Redistancing by flow time dependent {E}ikonal equation.
\newblock {\em J. Comput. Phys.}, 227(2):4002--4017, 2008.

\bibitem{chern_shu07}
I.-L. Chern and Y.-C. Shu.
\newblock A coupling interface method for elliptic interface problems.
\newblock {\em J. of Comput. Physics}, 225:2138--2174, 2007.

\bibitem{ciarlet78}
P.~G. Ciarlet.
\newblock {\em The Finite Element Method for Elliptic Problems}.
\newblock SIAM, 1978.

\bibitem{conti_niethammer_otto06}
S.~Conti, B.~Niethammer, and F.~Otto.
\newblock Coarsening rates in off-critical mixtures.
\newblock {\em SIAM J. Math. Anal.}, 37(6):1732--1741, 2006.

\bibitem{dillencourt_samet_tamminen92}
M.~Dillencourt, H.~Samet, and M.~Tamminen.
\newblock A general approach to connected-component labeling for arbitrary
  image representations.
\newblock {\em J. ACM}, 39, 1992.

\bibitem{dolbow_harari09}
J.~Dolbow and I.~Harari.
\newblock An efficient finite element method for embedded interface problems.
\newblock {\em J. Numer. Methods Eng.}, 78:229--252, 2009.

\bibitem{federer69}
H.~Federer.
\newblock {\em Geometric Measure Theory}.
\newblock Springer-Verlag, 1969.

\bibitem{folland76}
G.~B. Folland.
\newblock {\em Introduction to Partial Differential Equations}.
\newblock Princeton University Press, 1976.

\bibitem{gibou_fedkiw05}
F.~Gibou and R.~Fedkiw.
\newblock A fourth order accurate discretization for the laplace and heat
  equations on arbitrary domains, with applications to the stefan problem.
\newblock {\em J. Comput. Phys.}, 202:577--601, 2005.

\bibitem{gibou_fedkiw_cheng_kang02}
F.~Gibou, R.~Fedkiw, L.~Cheng, and M.~Kang.
\newblock A second order accurate symmetric discretization of the poisson
  equation on irregular domains.
\newblock {\em J. Comput. Phys.}, 176:1--23, 2002.

\bibitem{Giga:2006}
Y.~Giga.
\newblock {\em Surface evolution equations: A level set approach}, volume~99 of
  {\em Monographs in Mathematics}.
\newblock Birkh\"auser Verlag, Basel, 2006.

\bibitem{greenbaum_greengard_macfadden93}
A.~Greenbaum, L.~Greengard, and G.~B. MFadden.
\newblock Laplace's equation and the {D}irichlet-{N}eumann map in multiply
  connected domains.
\newblock {\em J. Comput. Phys.}, 105:267--278, 1993.

\bibitem{greengard_rokhlin87}
L.~Greengard and V.~Rokhlin.
\newblock A fast algorithm for particle simulations.
\newblock {\em J. Comput. Phys.}, 73(2):325--348, 1987.

\bibitem{gurtin86}
M.~Gurtin.
\newblock On the two-phase stefan problem with interfacial energy and entropy.
\newblock {\em Archive for Rational Mechanics and Analysis}, 96, 1986.

\bibitem{hansbo_hansbo02}
A.~Hansbo and P.~Hansbo.
\newblock An unfitted element method, based on nitsche's method for elliptic
  interface problems.
\newblock {\em Comput. Methods Appl. Mech. Eng.}, 191:5537--5552, 2002.

\bibitem{hansbo_hansbo04}
A.~Hansbo and P.~Hansbo.
\newblock A finite element method for the simulation of strong and weak
  discontinuities in solid mechanics.
\newblock {\em Comput. Methods Appl. Mech. Eng.}, 193:3523--3540, 2004.

\bibitem{huang_zou02}
J.~Huang and J.~Zou.
\newblock A mortar element method for elliptic problems with discontinous
  coefficients.
\newblock {\em IMA J. Numer. Anal.}, 22:549--576, 2002.

\bibitem{johansen97}
H.~Johansen.
\newblock {\em Cartesian grid embedded boundary finite difference methods for
  elliptic and parabolic differential equations on irregular domains}.
\newblock PhD thesis, University of California, Berkeley, 1997.

\bibitem{johansen_colella98}
H.~Johansen and P.~Colella.
\newblock A cartesian grid embedded boundary method for poisson's equation on
  irregular domains.
\newblock {\em J. Comput. Phys.}, 147:60--85, 1998.

\bibitem{karali_kevrekidis09}
G.~Karali and P.~Kevrekidis.
\newblock Bubble interactions for the {M}ullins-{S}ekerka problem: some case
  examples.
\newblock {\em Math. Comput. Simul.}, 80, 2009.

\bibitem{kress99}
R.~Kress.
\newblock {\em Linear Integral Equations}.
\newblock Springer-Verlag, New York, second edition, 1999.

\bibitem{kublik_tanushev_tsai13}
C.~Kublik, N.~M. Tanushev, and R.~Tsai.
\newblock An {I}mplicit {I}nterface {B}oundary {I}ntegral {M}ethod for
  {P}oisson's {E}quation on {A}rbitrary {D}omains.
\newblock {\em J. Comput. Phys.}, 247:269--311, 2013.

\bibitem{kublik_tsai16}
C.~Kublik and R.~Tsai.
\newblock Integration over curves and surfaces defined by the closest point
  mapping.
\newblock {\em Research in the Mathematical Sciences}, 2016.

\bibitem{leveque_li94}
R.~Leveque and Z.~Li.
\newblock The immersed interface method for elliptic equations with
  discontinuous coefficients and singular sources.
\newblock {\em SIAM J. Numer. Anal.}, 31:1019--1044, 1994.

\bibitem{li_ito06}
Z.~Li and K.~Ito.
\newblock The immersed interface method: numerical solutions of pdes involving
  interfaces and irregular domains (frontiers in applied mathematics).
\newblock {\em Society for Industrial and Applied Mathematics}, 2006.

\bibitem{liu_fedkiw_kang00}
X.~Liu, R.~Fedkiw, and M.~Kang.
\newblock A boundary condition capturing method for poisson's equation on
  irregular domains.
\newblock {\em J. Comput. Phys.}, 160(1):151--178, 2000.

\bibitem{mikhlin57}
S.~G. Mikhlin.
\newblock {\em Integral Equations}.
\newblock Pergammon, London, 1957.

\bibitem{Miura-2015}
T.-H. Miura.
\newblock Zero width limit of the heat equation on moving thin domains.
\newblock {\em UTMS Preprint Series, The University of Tokyo}, 2015.

\bibitem{mullins_sekerka63}
W.~W. Mullins and R.~F. Sekerka.
\newblock Morphological stability of a particle growing by diffusion or heat
  flow.
\newblock {\em J. of Applied Physics}, 34:323--329, 1963.

\bibitem{nystrom30}
E.~Nystr{\"{o}}m.
\newblock {\"{U}}ber die praktische {A}ufl{\"{o}}sung von {I}ntegralgleichungen
  mit {A}nwendungen auf {R}andwertaufgaben.
\newblock {\em Acta Math.}, 54:185--204, 1930.

\bibitem{osher_sethian88}
S.~Osher and J.~A. Sethian.
\newblock Fronts propagating with curvature dependent speed: Algorithms based
  on hamilton-jacobi formulations.
\newblock {\em J. Comp. Phys.}, 79:12--49, 1988.

\bibitem{peng_merriman_osher_zhao_kang99}
D.~Peng, B.~Merriman, S.~Osher, H.-K. Zhao, and M.~Kang.
\newblock A pde-based fast local level set method.
\newblock {\em J. Comput. Phys.}, 155(2):410--438, 1999.

\bibitem{russo_smereka00}
G.~Russo and P.~Smereka.
\newblock A remark on computing distance functions.
\newblock {\em J. Comput. Phys.}, 163:51--67, 2000.

\bibitem{sethian96}
J.~Sethian.
\newblock A fast marching level set method for monotonically advancing fronts.
\newblock {\em Proceedings of the National Academy of Sciences},
  93(4):1591--1595, 1996.

\bibitem{tsai_cheng_osher_zhao03}
Y.-H. Tsai, L.~Cheng, S.~Osher, and H.-K. Zhao.
\newblock Fast sweeping methods for a class of hamilton-jacobi equations.
\newblock {\em SIAM Journal on Numerical Analysis}, 41(2):673--694, 2003.

\bibitem{tsitsiklis95}
J.~Tsitsiklis.
\newblock Efficient algorithms for globally optimal trajectories.
\newblock {\em IEEE Transactions on Automatic Control}, 40:1528--1538, 1995.

\bibitem{zhu_chen_hou96}
J.~Zhu, X.~Chen, and T.~Y. Hou.
\newblock An efficient boundary integral method for the {M}ullins-{S}ekerka
  problem.
\newblock {\em J. Comput. Phys.}, 127:246--267, 1996.

\end{thebibliography}

\end{document}